\documentclass{gtmon_a}
\pdfoutput=1

\proceedingstitle{Groups, homotopy and configuration spaces (Tokyo
  2005)}
\conferencestart{5 July 2005}
\conferenceend{11 July 2005}
\conferencename{Groups, homotopy and configuration spaces,
                in honour of Fred Cohen's 60th birthday}
\conferencelocation{University of Tokyo, Japan}

\editor{Norio Iwase}
\givenname{Norio}
\surname{Iwase}

\editor{Toshitake Kohno}
\givenname{Toshitake}
\surname{Kohno}

\editor{Ran Levi}
\givenname{Ran}
\surname{Levi}

\editor{Dai Tamaki}
\givenname{Dai}
\surname{Tamaki}

\editor{Jie Wu}
\givenname{Jie}
\surname{Wu}

\title{Quasi-polynomials and the Bethe Ansatz}

\author{E Mukhin}
\givenname{E}
\surname{Mukhin}
\address{Department of Mathematical Sciences\\Indiana University --
Purdue University Indianapolis\\\newline
402 North Blackford St\\Indianapolis, IN 46202-3216\\USA}
\email{mukhin@math.iupui.edu}
\urladdr{}

\author{A Varchenko}
\givenname{A}
\surname{Varchenko}
\address{Department of Mathematics\\University of North Carolina 
at Chapel Hill\\\newline
Chapel Hill, NC 27599-3250\\USA}
\email{anv@email.unc.edu}
\urladdr{}

\volumenumber{13}
\issuenumber{}
\publicationyear{2008}
\papernumber{18}
\startpage{385}
\endpage{420}

\doi{}
\MR{}
\Zbl{}


\received{18 February 2006}
\revised{4 September 2006}
\accepted{11 October 2006}
\published{19 March 2008}
\publishedonline{19 March 2008}
\proposed{}
\seconded{}
\corresponding{}
\version{}

\arxivreference{math.QA/0604048}

\makeatletter
\def\cnewtheorem#1[#2]#3{\newtheorem{#1}{#3}[section]
\expandafter\let\csname c@#1\endcsname\c@theorem}

\let\xysavmatrix\xymatrix
\def\xymatrix{\disablesubscriptcorrection\xysavmatrix}
\AtBeginDocument{\let\bar\wbar\let\tilde\wtilde\let\hat\what}
\makeop{XXX}

\newtheorem{theorem}{Theorem}[section]
\cnewtheorem{conj}[theorem]{Conjecture}
\cnewtheorem{cor}[theorem]{Corollary}
\cnewtheorem{lem}[theorem]{Lemma}
\cnewtheorem{prop}[theorem]{Proposition}

\theoremstyle{definition}

\theoremstyle{remark}
\cnewtheorem{rem}[theorem]{Remark}
\theoremstyle{remark}
\cnewtheorem{example}[theorem]{Example}

\makeatother  

\newcommand{\nc}{\newcommand}
\nc{\on}{\operatorname}
\nc{\ch}{\mbox{ch}}
\nc{\pone}{{\mathbb C}{\mathbb P}^1}
\nc{\pa}{\partial}
\nc{\F}{{\mathcal F}}
\nc{\arr}{\rightarrow}
\nc{\larr}{\longrightarrow}
\nc{\al}{\alpha}
\nc{\ri}{\rangle}
\nc{\lef}{\langle}
\nc{\W}{{\mathcal W}}
\nc{\la}{\lambda}
\nc{\ep}{\epsilon}

\nc{\Om}{\Omega}
\nc{\su}{\widehat{{\mathfrak sl}}_2}
\nc{\sw}{{\mathfrak s}{\mathfrak l}}

\nc{\g}{{\mathfrak g}}
\nc{\h}{{\mathfrak h}}
\nc{\n}{{\mathfrak n}}
\nc{\N}{\widehat{\n}}
\nc{\G}{\widehat{\g}}
\nc{\De}{\Delta_+}
\nc{\Ga}{\Gamma}
\nc{\one}{{\mathbf 1}}
\nc{\z}{{\mathfrak Z}}
\nc{\zz}{{\mathcal Z}}
\nc{\Hh}{{\mathcal H}_\beta}
\nc{\qp}{q^{\frac{k}{2}}}
\nc{\qm}{q^{-\frac{k}{2}}}
\nc{\La}{\Lambda}
\nc{\wt}{\widetilde}
\nc{\qn}{\frac{[m]_q^2}{[2m]_q}}
\nc{\cri}{_{\on{cr}}}
\nc{\kk}{h^\vee}
\nc{\sun}{\widehat{\sw}_N}
\nc{\hh}{\widehat{\mathfrak h}}
\nc{\HH}{{\mathcal H}_{q,t}}
\nc{\ca}{\wt{{\mathcal A}}_{h,k}(\sw_2)}
\nc{\gl}{\widehat{{\mathfrak g}{\mathfrak l}}_2}
\nc{\el}{\ell}
\nc{\s}{{\mathbf s}}
\nc{\bi}{\bibitem}
\nc{\om}{\omega}
\nc{\WW}{\W_\beta}
\nc{\scr}{{\mathbf S}}
\nc{\ab}{{\mathbf a}}
\nc{\rr}{r}
\nc{\ol}{\overline}
\nc{\con}{qt^{-1} + q^{-1}t}
\nc{\den}{q^{\el-1} t^{-\el+1}+ q^{-\el+1} t^{\el-1}}
\nc{\ds}{\displaystyle}
\nc{\B}{B}
\nc{\A}{{\mathbb A}}
\nc{\GG}{{\mathcal G}}
\nc{\UU}{{\mathcal U}}
\nc{\MM}{{\mathcal M}}
\nc{\CC}{{\mathcal C}}
\nc{\GL}{{}^L G}
\nc{\dzz}{\frac{dz}{z}}
\nc{\Res}{\on{Res}}
\nc{\rep}{{\mathcal R}ep \;}
\nc{\uqg}{U_q \G}
\nc{\uqgg}{U_q \g}
\nc{\Fq}{{\mathbb F}_q}

\nc{\stimes}{\ltimes}
\nc{\K}{\hat{\mathcal K}}
\nc{\Ql}{\ol{\mathbb Q}_\ell}

\nc{\ga}{\gamma}
\nc{\PL}{{}^L P}
\nc{\E}{\mc E}
\nc{\mc}{\mathcal}
\nc{\bb}{{\mathfrak b}}
\nc{\OO}{{\mc O}}
\nc{\Po}{{\mc P}}
\nc{\V}{{\mc V}}
\nc{\yy}{{\mc Y}}
\nc{\M}{\mathcal M}
\nc{\Coh}{{{\mathcal C}oh}}
\nc{\Cohn}{\Coh_n}
\nc{\f}{{\mathcal F}}
\nc{\si}{_E}
\nc{\Gaf}{{\mathbb G}_{a,\Fq}}
\nc{\KK}{{\mathfrak k}}

\nc{\PCr}{{ \bs P  (\C[x])^r   }}
\nc{\PCN}{{ \bs P  (\C[x])^N   }}
 
\nc{\sN}{sl_{2N+1}}

\nc{\Pzr}{{ \bs P(\C((x-z)))^r}}

\nc{\PzN}{{ \bs P(\C((x-z)))^N}}

\newcommand{\bean}{\begin{eqnarray}}
\newcommand{\eean}{\end{eqnarray}}
\newcommand{\be}{\begin{displaymath}}
\newcommand{\ee}{\end{displaymath}}
\newcommand{\bea}{\begin{eqnarray*}}   
\newcommand{\eea}{\end{eqnarray*}}
\newcommand{\bs}{\boldsymbol}

\begin{document}

\begin{asciiabstract}
  We study solutions of the Bethe Ansatz equation related to the
  trigonometric Gaudin model associated to a simple Lie algebra g
  and a tensor product of irreducible finite-dimensional
  representations.  Having one solution, we describe a construction of
  new solutions.  The collection of all solutions obtained from a
  given one is called a population.  We show that the Weyl group of
  g acts on the points of a population freely and transitively (under
  certain conditions).
  
  To a solution of the Bethe Ansatz equation, one assigns a common
  eigenvector (called the Bethe vector) of the trigonometric Gaudin
  operators.  The dynamical Weyl group projectively acts on the common
  eigenvectors of the trigonometric Gaudin operators. We conjecture
  that this action preserves the set of Bethe vectors and
  coincides with the action induced by the action on points of
  populations. We prove the conjecture for sl_2.
\end{asciiabstract}

\begin{htmlabstract}
<p class="noindent">
We study solutions of the Bethe Ansatz equation related to the
trigonometric Gaudin model associated to a simple Lie algebra g
and a tensor product of irreducible finite-dimensional
representations.  Having one solution, we describe a construction of
new solutions.  The collection of all solutions obtained from a
given one is called a population.  We show that the Weyl group of
g acts on the points of a population freely and transitively (under
certain conditions).
</p>
<p class="noindent">
To a solution of the Bethe Ansatz equation, one assigns a common
eigenvector (called the Bethe vector) of the trigonometric Gaudin
operators.  The dynamical Weyl group projectively acts on the common
eigenvectors of the trigonometric Gaudin operators. We conjecture
that this action preserves the set of Bethe vectors and
coincides with the action induced by the action on points of
populations. We prove the conjecture for sl<sub>2</sub>.
</p>
\end{htmlabstract}

\begin{abstract}
  We study solutions of the Bethe Ansatz equation related to the
  trigonometric Gaudin model associated to a simple Lie algebra $\mathfrak{g}$
  and a tensor product of irreducible finite-dimensional
  representations.  Having one solution, we describe a construction of
  new solutions.  The collection of all solutions obtained from a
  given one is called a population.  We show that the Weyl group of
  $\mathfrak{g}$ acts on the points of a population freely and transitively (under
  certain conditions).
  
  To a solution of the Bethe Ansatz equation, one assigns a common
  eigenvector (called the Bethe vector) of the trigonometric Gaudin
  operators.  The dynamical Weyl group projectively acts on the common
  eigenvectors of the trigonometric Gaudin operators. We conjecture
  that this action preserves the set of Bethe vectors and
  coincides with the action induced by the action on points of
  populations. We prove the conjecture for $sl_2$.
\end{abstract}

\maketitle

\section{Introduction}
The Bethe Ansatz is a method to diagonalize a commuting family of
linear operators, usually called Hamiltonians.  The method is applied
to Hamiltonians of numerous quantum integrable systems.  Given a
solution of a suitable system of equations (called the Bethe Ansatz
equation), the Bethe Ansatz produces an eigenvector (called the Bethe
vector).  This paper is motivated by the Bethe Ansatz method applied
to the trigonometric Gaudin model, see Markov, Schechtman and the
second author \cite{MaV,SV}, and \fullref{act sec}.

For the case of the trigonometric Gaudin model the Bethe Ansatz equation
and the Bethe vectors depend on an additional parameter, a generic
$\g$--weight $\la$.  The Bethe Ansatz equation has the form \eqref{Bethe
  eqn}. The Bethe vectors have the form \eqref{bethe vector}, see
\fullref{Bethe vector prop}.  The Bethe Ansatz equation
\eqref{Bethe eqn} can be formulated as a system of suitable Wronskian
equations for a tuple of polynomials $\bs y=(y_1,\dots,y_r)$ of one
variable, where $r$ is the rank of $\g$ and the polynomials are
labeled by simple roots of $\g$, see \fullref{gen rep}.

For example, let $\g = sl_2$. The $sl_2$--weights can be
identified with complex numbers. Consider the
trigonometric Gaudin model associated to the tensor product of
irreducible $sl_2$--modules with highest weights $\La_j$, located 
respectively at points
$z_j$. In this case,
the Bethe Ansatz equation with parameter $\la \in \C$ is an equation
on one polynomial $y$. The polynomial $y$ satisfies the Bethe Ansatz
equation, if and only if its roots are simple and there exists another
polynomial $\tilde y$ such that
\be 
y'(x^{\la+1}\tilde y)\ -\ y(x^{\la+1}\tilde
y)'\ =\ x^\la\prod_j(x-z_j)^{\La_j}\ .  
\ee 
For a given $y$ and a non-integer $\la$, the polynomial $\tilde y$
is unique.  One can show that for
 almost all $\la$, the roots of $\tilde y$ are simple.
Moreover, if the roots of $\tilde y$ are simple, then the
polynomial $\tilde y$ also satisfy the Bethe Ansatz equation but
with the
new parameter $-\la-2$.  Thus from one solution of the Bethe Ansatz
equation with parameter $\la$ (the polynomial $y$) we obtain another
solution with parameter $-\la-2$ (the polynomial $\tilde y$).  We call this
procedure the simple reproduction procedure.

For an arbitrary simple Lie algebra $\g$, there is a similar simple
reproduction procedure associated with every simple root of $\g$. 
Consider an $r$--tuple $\bs y=(y_1,\dots,y_r)$ of polynomials forming
 a solution of the Bethe Ansatz
equation associated with a generic $\g$--weight $\la$.  Then we have
the $i$-th simple reproduction procedure for $i=1,\dots,r$. The
$i$-th simple reproduction procedure constructs a new tuple
$\bs y^{(i)}=(y_1,\dots,y_{i-1}, \tilde y_i, y_{i+1},\dots, y_r)$
under certain conditions.

We call an $r$--tuple of polynomials $\bs y=(y_1,\dots,y_r)$ fertile
with respect to $\la$ if the $i$-th simple reproduction procedure is
well-defined for $i=1,\dots,r$.  In particular, if $\bs y$ forms a
solution of the Bethe Ansatz equation associated to $\la$, then $\bs
y$ is fertile with respect to $\la$.  Moreover, if the $i$-th simple
reproduction procedure results in a generic (in an appropriate sense)
$r$--tuple $\bs y^{(i)}$, then $\bs y^{(i)}$ also forms a solution of
the Bethe Ansatz equation associated to the weight $s_i\cdot\la$, where
$s_i$ is the $i$-th elementary reflection in the Weyl group of $\g$. 
It follows that the $r$--tuple $\bs y^{(i)}$ is fertile with respect to
$s_i\cdot\la$.

We call an $r$--tuple of polynomials $\bs y$ super-fertile with respect
to $\la$ if all iterations of the simple reproduction procedures are well
defined.  We conjecture that if $\bs y$ forms a solution of
the Bethe Ansatz equation then $\bs y$ is super-fertile.  We prove
the conjecture for simple Lie algebras of type $A_r,B_r$.

The set of all $r$--tuples obtained from a given super-fertile
$r$--tuple by iterations of simple reproduction procedures is called a
population.

For simple Lie algebras, we prove that the population obtained from a
super-fertile $r$--tuple associated to a generic weight $\la$ contains
exactly one $r$--tuple associated to every weight of the form $w\cdot
\la$, where $w$ runs through the elements of the Weyl group of $\g$.
We also prove that the population does not contain any other
$r$--tuples. This one-to-one correspondence between the tuples of the
population and the weights of the form $w\cdot\la$, 
allows us to introduce a free and
transitive action of the Weyl group on points of the population. Then
the action of simple reflections is given by the simple reproduction
procedures. The proof is based on the important fact that in the case
of $A_r$, the populations are in one-to-one
correspondence with certain spaces of quasi-polynomials, see \fullref{sing points}.

If all elements of a population are generic and therefore correspond
to solutions of the Bethe Ansatz equation, then we have an action of
the Weyl group on the set of the solutions.  In particular this
defines an action of the Weyl group on the set of the associated Bethe
vectors, considered up to proportionality.

On the other hand, the dynamical Weyl group commutes with the
trigonometric Gaudin operators and projectively acts on eigenvectors
of the trigonometric Gaudin operators, see Tarasov and the second
author \cite{TV} and \fullref{proper}.

We conjecture that the action of the dynamical Weyl group maps the
Bethe vectors to (scalar multiples of) the Bethe vectors and moreover,
the two actions on the Bethe vectors coincide. We prove this
conjecture for $\g = sl_2$.

The reproduction procedure exists for solutions of the Bethe Ansatz
equations associated with many quantum integrable models.  In Sections
\ref{exp sec}, \ref{xxx sec} we give two other examples of the
situation in which the corresponding Bethe Ansatz equation admits a
reproduction procedure and prove that the elements of the
corresponding populations are also labeled by the elements of the Weyl
group. These examples are related to the quasi-periodic Gaudin and
$\XXX$ models.  In joint work with Tarasov \cite{MTV2}, we apply the
results of Sections \ref{exp sec} and \ref{xxx sec} to study the Bethe
Ansatz of $(gl_N,gl_M)$--dual quasiperiodic Gaudin and $\XXX$ models.

The notions of the reproduction procedure and populations for the
Bethe Ansatz equation of the Gaudin and $\XXX$--type models were
introduced in our papers \cite{MV,MV2}, see also \cite{MV3,MV5}.  The
populations in that situation are shown to be isomorphic to the flag
variety of the Langlands dual algebra $\g^\vee$ for all simple Lie
algebras in the case of the Gaudin model by Frenkel \cite{F} and the
authors \cite{MV4}, and for the $\XXX$--type model by the authors
\cite{MV5}.

The paper is constructed as follows. \fullref{crit point sec} contains
notation and definitions.  In \fullref{pop sec} we define the
reproduction procedure and populations. In \fullref{proof sec} we
prove that for simple Lie algebras, the Weyl group acts freely and
transitively on the elements of a population.  In \fullref{act sec} we
discuss two actions of the Weyl group on the Bethe vectors, the one
given by the action of the dynamical Weyl groups and the one given by
the reproduction procedure. We compare them in the case of $sl_2$. In
Sections \ref{exp sec} and \ref{xxx sec} we describe two more examples
of the situation in which the Bethe Ansatz admit a reproduction
procedure and the the Weyl group acts freely and transitively on the
elements of a population.

\medskip
{\bf Acknowledgements}\qua
We thank V Tarasov for many valuable discussions.

Research of EM is supported in part by NSF grant DMS-0601005.
Research of AV is supported in part by NSF grant DMS-0244579.

\section{Master functions and critical points}\label{crit point sec}
\subsection{Kac--Moody algebras}\label{Kac_Moody sec}
Let $A=(a_{ij})_{i,j=1}^r$ be a generalized Cartan matrix, $a_{ii}=2$,
$a_{ij}=0$ if and only $a_{ji}=0$, $a_{ij}\in \Z_{\leq 0}$ if $i\ne
j$.  We assume that $A$ is symmetrizable, there is a diagonal matrix
$D=\on{diag}\{d_1,\dots,d_r\}$ with positive integers $d_i$ such that
$B=DA$ is symmetric.

Let $\g=\g(A)$ be the corresponding complex Kac--Moody Lie algebra (see
Kac \cite[Section 1.2]{K}), $\h \subset \g$ the Cartan subalgebra.  The
associated scalar product is non-degenerate on $\h^*$ and $\on{dim}\h
= r + 2d$, where $d$ is the dimension of the kernel of the Cartan
matrix $A$.

Let $\al_i\in \h^*$, $\al_i^\vee\in \h$, $i = 1, \dots , r$, be the
sets of simple roots, coroots, respectively. We have
\bea
 (\al_i,\al_j)&=& d_i \ a_{ij}, \\
\langle\la ,\al^\vee_i\rangle&=
&2(\la,\al_i)/{(\al_i,\al_i)}, \qquad \la\in\h^*.
\eea
Let $\mathcal P = \{ \lambda \in \h^* \, |\, \langle\la
,\al^\vee_i\rangle \in \Z\}$ and $\mathcal P^+ = \{ \lambda \in \h^*
\, |\, \langle\la ,\al^\vee_i\rangle \in \Z_{\geq 0}\}$ be the sets of
integral and dominant integral weights.

Fix $\rho\in\h^*$ such that $\langle\rho,\al_i^\vee\rangle=1$,
$i=1,\dots,r$. We have $(\rho,\al_i)= (\al_i,\al_i)/2$.

The Weyl group $\mathcal W\in\on{End (\h^*)}$ is generated by 
reflections $s_i$, $i=1,\dots,r$, 
\be
s_i(\la)=\la-\langle\la,\al_i^\vee\rangle\al_i, \qquad \la\in\h^*.
\ee
We use the notation
\bea\label{shifted}
w\cdot\la=w(\la+\rho)-\rho,\qquad w\in \mathcal W,\;\la\in\h^*,
\eea
for the shifted action of the Weyl group.

\subsection{The definition of
  master functions and critical points}\label{master sec} 
We fix a Kac--Moody algebra $\g=\g(A)$. We fix $\bs\La=(\La_i)_{i=1}^n$,
$\La_i\in\mathcal P^+$; $\bs z=(z_i)_{i=1}^n\in\C^n$. We assume
$z_i\neq 0$ and $z_i\neq z_j$ if $i\neq j$. The parameters $\bs \La,
\bs z$ are always fixed and we often do not stress the dependence of
our objects on these parameters.

In addition we choose $\la\in\h^*$ and $\bs
l=(l_1,\dots,l_r)\in\Z^r_{\geq 0}$.  The choice of $\bs l$ is
equivalent to the choice of the weight at infinity $\La_\infty$
defined by the formula:
\bean\label{at inf}
\La_\infty\ =\sum_{s=1}^n \La_s  -  \sum_{i=1}^r
l_i\al_i \in \mathcal P .
\eean
The {\it master function} $\Phi(\bs t;\La_\infty,\la)$
is defined as follows (see Felder, Schechtman and the second author
\cite{FV,SV}):
\begin{gather}\label{master}
\Phi(\bs t; \La_\infty,\la)= 
\prod_{i=1}^r\prod_{j=1}^{l_i}(t_j^{(i)})^{-(\la,\al_i)}\
\prod_{i=1}^r\prod_{j=1}^{l_i}\prod_{s=1}^n
  (t_j^{(i)}-z_s)^{-(\La_s,\al_i)}\ \times\\
\prod_{i=1}^r\prod_{1\leq j<s\leq l_i} (t_j^{(i)}-t_s^{(i)})^{(\al_i,\al_i)}
  \prod_{1\leq i<j\leq r}
\prod_{s=1}^{l_i}\prod_{k=1}^{l_{j}}(t_s^{(i)}-t_k^{(j)})^{(\al_i,\al_j)}.
\notag
\end{gather}
The master function $\Phi$ is a function of variables 
$\bs  t = ( t_j^{(i)})^{j = 1, \dots , l_i}_{i = 1,\dots,r}$.

The master function $\Phi$ is symmetric with respect to permutations
of variables with the same upper index.

A point $\bs t$ with complex coordinates is called a {\it critical
  point} associated to $(\bs \La,\bs z;\La_\infty$, $\la)$ (we often
write just $(\La_\infty, \la)$) if the following system of algebraic
equations is satisfied
\bean\label{Bethe eqn}
-\frac{(\la,\al_i)}{t_j^{(i)}}-\sum_{s=1}^n 
\frac{(\Lambda_s, \alpha_i)}{t_j^{(i)}-z_s}\ +\
\sum_{s,\ s\neq i}\sum_{k=1}^{l_s} \frac{(\alpha_s, 
\alpha_i)}{ t_j^{(i)} -t_k^{(s)}}\ +\
\sum_{s,\ s\neq j}\frac {(\alpha_i, \alpha_i)}{ t_j^{(i)} -t_s^{(i)}}
= 0, 
\eean
where $i = 1, \dots , r$, $j = 1, \dots , l_i$.
In other words, a point $\bs t$ is a critical point if
\be
\left(\Phi^{-1}\frac{\partial \Phi }{\partial t_j^{(i)}}\right)
(\bs t;\La_\infty,\la)=0, \qquad
i=1,\dots,r,\; j=1,\dots, l_i.
\ee
Note that the product of symmetric groups $S_{\bs l}=S_{l_1}\times
\dots \times S_{l_r}$ acts on the critical set of the master function
permuting the coordinates with the same upper index. All orbits have
the same cardinality $l_1! \dots l_r!$.  We make no distinction
between critical points in the same orbit.

The system of equations \eqref{Bethe eqn} coincides with the Bethe
Ansatz equation of the non-homogeneous Gaudin model, see Reshetikhin
and the second author \cite{RV} and \fullref{Bethe vector prop}.

\begin{lem}\label{finite}
  For almost all values of $\la$, the set of critical points
  associated to $(\La_\infty,\la)$ is finite.
\end{lem}
\begin{proof}
  The lemma follows from \cite[Lemma 2.1]{MV}.
\end{proof}
 
\section{Populations}\label{pop sec}

\subsection{The $sl_2$ populations}
In the case of $sl_2$ the dominant integral weights are identified
with non-negative integers and the system of equations \eqref{Bethe eqn}
takes the form:
\bean \label{bethe sl2}
-\frac{\la}{t_j}-\sum_{s=1}^n\frac{\La_s}{t_j-z_s}+\sum_{k,\ k\neq
  j}^l\frac{2}{t_j-t_k}=0,
\eean
$j=1,\dots, l$, where $\La_s\in\Z_{>0}$ and $\la\in\C$. The weight $\La_\infty$ is given by $\La_\infty=\sum_{s=1}^n\La_s-2l$.

We set
\be
T(x)=\prod_{s=1}^n (x-z_s)^{\La_s}, \qquad F=x\prod_{s=1}^n (x-z_s).
\ee
Given a tuple $\bs t=(t_1,\dots,t_l)$ we represent it by the polynomial
$y(x)=\prod_{j=1}^l (x-t_j)$. We are interested in the zeros of $y(x)$
and therefore we make no distinction between $y(x)$ and $cy(x)$,
where $c$ is a non-zero constant.

A polynomial $y(x)$ is called {\it off-diagonal} with respect to $(\bs
\La,\bs z)$ if $y(x)$ has only simple roots, $y(0)\neq 0$ and
$y(z_s)\neq 0$ for all $s = 1, \dots , n$ such that $\La_s\neq 0$.

Since $(\bs \La,\bs z)$ is fixed, we often call polynomials $y(x)$
off-diagonal with respect to $(\bs \La,\bs z)$ simply off-diagonal.

\begin{lem}{\rm(T\,J Stieltjes \cite[Section 6.81]{S})}\label{st}\qua
  A polynomial $y$ of degree $l$ represents an $sl_2$ critical point
  associated to $(\La_\infty,\la)$ if and only if $y$ is off-diagonal
  and there exists a polynomial $C(x)$ such that
\bean\label{de}
F(x)\,y''- F(x)\,\ln' (x^\la T(x))\,
y'(x)+C(x)\,y(x)=0. 
\eean
\end{lem}
\begin{proof}
  Equation \eqref{bethe sl2} can be reformulated as the statement that
  the function $y''-\ln' (x^\la T)y'$ equals zero at $x=t_j$ for all
  $j$. Therefore \eqref{bethe sl2} is equivalent to the divisibility of
  the polynomial $Fy''-F\ln' (x^\la T)y'$ by $y(x)$.
\end{proof}

Note that the coefficients of $y'',y'$, and $y$ in \eqref{de} are
polynomials of degree $n+1,n$, and at most $n-1$, respectively.

\begin{theorem}\label{sl2 rep} Let  
  $y$ represent an $sl_2$ critical point associated to $(\bs \La,\bs
  z;\La_\infty$, $\la)$. Then equation \eqref{de} has a solution of the
  form $x^{\la+1} \tilde y(x)$ where $\tilde y(x)$ is a polynomial. If
  $\la$ is not a negative integer, then such a polynomial $\tilde
  y(x)$ is unique.
  
  Moreover, there exists a finite set $C(\bs\La,\bs
  z;\La_\infty)\subset \C$ such that if $\la$ is not a negative
  integer, $\la\not\in C(\bs \La,\bs z;\La_\infty)$, and $y$
  represents an $sl_2$ critical point associated to $(\bs \La,\bs
  z;\La_\infty$, $\la)$ then polynomial $\tilde y(x)$ represents an
  $sl_2$ critical point associated to $(\bs \La,\bs z; -\La_\infty$,
  $-\la-2)$.
\end{theorem}

We call $\tilde y(x)$ {\it the immediate descendent of $y(x)$ with
  respect to $\lambda$.}

\begin{proof}
  Equation \eqref{de} is a Fuchsian differential equation, with singular
  points at $0, z_1,\dots,$ $z_n,$ $ \infty$.  At $0$ the exponents of
  the equation are $0,\la+1$. Therefore, around $0$ there is a
  solution of the form $u(x)=x^{\la+1} \tilde y(x)$ where $\tilde
  y(x)$ is a function holomorphic and non-vanishing at $x=0$. Such a
  solution is unique if $\la+1\not\in\Z_{\le 0}$.
  
  At a point $z_s$ the exponents are $(0,\La_s+1)$.  Since $y(z_s)\neq
  0$ and $y$ is a polynomial solution, there is no monodromy around
  $z_s$, and thus $\tilde y(x)$ is an entire function, cf Scherbak
  \cite[Lemma 7]{ScV}.  The function $\tilde y(x)$ is a polynomial
  since equation \eqref{de} is Fuchsian.
\be
W(f,g)=f'g-fg'\tag*{\hbox{Denote}}
\ee
the Wronskian of functions $f$ and $g$.
We have 
\bean\label{aux}
W(y,x^{\la+1} \tilde y(x))=x^\la T(x).
\eean
\be
\deg y+\deg \tilde y \ = \ \sum_{s=1}^n\La_s \ .\tag*{\hbox{Thus}}
\ee
The polynomial $\tilde y(x)$ satisfies the equation
\be
F(x)\,\tilde y''- F(x)\, \ln' (x^{-\la-2}T(x))\,
\tilde y'(x)+\tilde C(x)\,\tilde y(x)=0. 
\ee
Thus if
$\tilde y(x)$ is off-diagonal, then the polynomial
$\tilde y(x)$  represents a critical point
associated to $(\bs \La, \bs z; -\La_\infty,-\la-2)$.

Finally, we prove that for all but finitely many $\la$, the polynomial
$\tilde y(x)$ is off-diagonal. If $\tilde y(x)$ is not off-diagonal,
then $\tilde y(x)$ has zero of order $\La_s+1$ at at least one of $z_s$.
We show that such a pair $(y(x), \tilde y(x))$ is possible for at most
finitely many $\la$.

Consider a family of polynomials $y_\la(x)=\prod_j^l(x-t_{j,\la})$
which algebraically depends on $\la$ and such that $y_\la(x)$
represents a critical point corresponding to $(\La_\infty,\la)$ for
all but finitely many $\la$.
 
We have a finite number of such families and for all but finitely many
$\la$ every polynomial representing a critical point belongs to such a
family.

Let $\tilde y_\la(x)=\prod_{j=1}^{\tilde l}(x-\tilde t_{j,\la})$ be the 
descendent polynomial of $y_\la(x)$. 

From \eqref{aux} we get
\be
\frac{\la+1}{x}+\sum_{i=1}^{\tilde l} \frac{1}{x-\tilde
  t_{j,\la}}-\sum_{i=1}^l
\frac{1}{x-t_{j,\la}}=\frac{(\la+\tilde l-l)\,\prod_{j=1}^n(x-z_j)^{\La_j}}
{x\prod_{i=1}^{\tilde l} (x-\tilde t_{j,\la}) \ \prod_{i=1}^l(x-t_{j,\la})}.
\ee
Let $\la$ tend to infinity.  Comparing the main terms of asymptotics
of the left and right hand sides , we conclude that the limit of
$y_\lambda\tilde y_\lambda$ is $T$.  The polynomial $T$ has zero of
order $\La_s$ at $z_s$ and therefore $\tilde y_\la$ cannot have a zero
of order $\La_s+1$ at $z_s$ for all but finitely many $\la$.  It
follows that $\tilde y_\la$ has a zero at $z_s$ for only finitely many
values of $\la$.
\end{proof}

\begin{cor}\label{sl2 fertile}
  A polynomial $y(x)$ represents an $sl_2$ critical point associated
  to $(\La_\infty,\la)$ if and only if $\deg y(x)=l$, $y(x)$ is
  off-diagonal with respect to $(\bs \La,\bs z)$ and there exists a
  polynomial $\tilde y(x)$ such that $W(y,x^{\la+1}\tilde y)=x^{\la}
  T(x)$.
\end{cor}

Note that if $y(x)$ represents a critical point associated to
$(\La_\infty,\la)$ and if the descendent polynomial $\tilde y(x)$ is
off-diagonal then $\tilde y(x)$ represents a critical point associated
to $(w\La_\infty, w\cdot \la)$, where $w\neq \on{id}$ is the generator of the
$sl_2$ Weyl group.

The polynomials $y(x)$ and $\tilde {y}(x)$ may coincide (up to a
multiplicative constant). For example, if there are no $z_s$, if $l=0$
and $y(x)=1$, then $\tilde y(x) = 1/(\la+1)$. The polynomials $y(x)$
and $\tilde y(x)$ are constant multiples of each other, but they
represent critical points associated to different weights
$(\La_\infty=0,\,\la)$ and $(\La_\infty=0,\,-\la-2)$, respectively.

Assume that $\la\in\C$ is not an integer.  Then the unordered pair
$\{\,(y(x),\la),\, $ $(\tilde {y}(x),-\la-2)\,\}$, is called {\it the
  $sl_2$ population originated at $(y(x),\la)$}.  Here $y(x),\tilde
y(x)$ are considered up to a multiplicative constant.

\begin{lem}\label{ss} Let $\la$ be not an integer.
  Let $P=\{\,(y(x),\la),\, (\tilde y(x),-\la-2)\,\}$ be the $sl_2$
  population originated at $(y(x),\la)$.  Assume that $\tilde y(x)$ is
  off-diagonal.  Then $y(x)$ is a descendent polynomial of $\tilde
  y(x)$ with respect to $-\la-2$ and the population originated at
  $(\tilde y(x),-\la-2)$ coincides with $P$.
\end{lem}
\begin{proof} We have
\be
x^{-\la-2}T=(x^{-\la-1})^2x^\la T=(x^{-\la-1})^2W(y,x^{\la+1}\tilde
y)=W(\tilde y,-x^{-\la-1}y),
\ee
cf the Wronskian identities in \cite[Lemma 9.2]{MV}.
\end{proof}

\subsection{The populations in the case of a general Kac--Moody algebra}
\label{gen pop sec}
Recall that we fixed a Kac--Moody algebra $\g=\g(A)$ of rank $r$, a
non-negative integer $n\in\Z_{\geq 0}$, an $n$--tuple of $\g$--weights
$\bs\La=(\La_i)_{i=1}^n$, $\La_i\in\mathcal P^+$, and an $n$--tuple
$\bs z=(z_i)_{i=1}^n\in\C^n$ such that $z_i\neq 0$ and $z_i\neq z_j$
if $i\neq j$.

For $i=1,\dots,r$, we set 
\bean\label{T}
T_i(x)=\prod_{s=1}^n
(x-z_s)^{\langle \La_s,\al_i^\vee\rangle}.
\eean
Given a set of numbers $\bs
t=(t_j^{(i)})_{i=1,\dots,r}^{j=1,\dots,l_i}$, we represent it by the
$r$--tuple of polynomials $\bs y=(y_1,\dots,y_r)$, where
$y_i(x)=\prod_{j=1}^{l_i}(x-t_j^{(i)})$, $i=1,\dots,r$.  We are
interested only in the roots of the polynomials $y_i$, therefore we
make no distinction between tuples $(y_1,\dots,y_r)$ and
$(c_1y_1,\dots,c_ry_r)$, where $c_i$ are non-zero constants.

An $r$--tuple of polynomials $\bs y$ is called {\it off-diagonal with
  respect to $(\bs \La, \bs z)$} if its roots do not belong to the
union of singular hyperplanes in \eqref{Bethe eqn}.  Namely $\bs y$ is
off-diagonal if for $i=1,\dots,n$, all roots of polynomial $y_i$ are
simple, non-zero, different from the roots of polynomials $y_j$ for
all $j$ such that $(\al_j,\al_i)\neq 0$ and different from the roots
of polynomial $T_i$.

Since $\bs \La, \bs z$ are fixed, we often call $r$--tuples of
polynomials $\bs y$ which are off-diagonal with respect to $(\bs \La,
\bs z)$ simply off-diagonal.

An $r$--tuple of polynomials $\bs y$ is called {\it fertile in the
  $i$-th direction, $i\in\{1,\dots,r\}$ with respect to $\la$} if
there exists a polynomial $\tilde y_i$ such that
\bean\label{wr}
W(y_i,x^{\langle  \la+\rho,\al_i^\vee\rangle }\tilde y_i)=x^{\langle
  \la,\al_i^\vee \rangle}\ T_i\prod_{j,\;j\neq i} y_j^{-a_{ij}}.
\eean
Then the tuple of polynomials $\bs y^{(i)}:=(y_1,\dots,\tilde
y_i,\dots,y_r)$ is called {\it an immediate descendent of $\bs y$ in
  the $i$-th direction with respect to $\la$}.

Recall that $s_i\in\mc W$ are elementary reflections in the Weyl group
of $\g$.

We call an $r$--tuple of polynomials $\bs y^{(i_1,i_2,\dots,i_m)}$,
where $i_k\in\{1,\dots,r\}$, $k=1,\dots,m$, {\it a descendent of $\bs
  y$ with respect to $\la$ in the directions $(i_1,\dots,i_m)$} if
there exist $r$--tuples of polynomials $\bs y^{(i_1,i_2,...,i_k)}$,
$k=1,\dots,m-1$, such that for $k=1,\dots,m$, the $r$--tuple $\bs
y^{(i_1,i_2,...,i_k)}$ is an immediate descendent of $\bs
y^{(i_1,i_2,...,i_{k-1})}$ in the $i_k$-th direction with respect to
$(s_{i_{k-1}}\dots s_{i_2} s_{i_1})\cdot \la$.

An $r$--tuple of polynomials $\bs y$ is called {\it fertile with
  respect to $\la$} if it is fertile in all directions $i=1,\dots,r$.

An $r$--tuple of polynomials $\bs y$ is called {\it super-fertile with
  respect to $\la$} if it is fertile with respect to $\la$ and all
descendents $\bs y^{(i_1,i_2,...,i_k)}$ of $\bs y$ with respect to
$\la$ are fertile with respect to $(s_{i_{k}}\dots s_{i_1})\cdot \la$.

\begin{theorem}\label{gen rep} Let $\g$ be a Kac--Moody algebra.
  An $r$--tuple $\bs y$ represents a $\g$ critical point associated to
  $(\bs \La,\bs z;\La_\infty, \la)$ if and only if $\deg y_i=l_i$,
  $i=1,\dots,r$, $\bs y$ is off-diagonal with respect to $(\bs\La,\bs
  z)$ and fertile with respect to $\la$.  Moreover, if an 
  immediate descendent of $\bs y$ in the direction $i$,
  $\bs y^{(i)}=(y_1,\dots,\tilde y_i,\dots,y_r)$, is off-diagonal with
  respect to $(\bs\La,\bs z)$ then it represents a $\g$ critical point
  associated to $(\bs \La,\bs z;s_i\La_\infty,s_i\cdot \la)$.
\end{theorem}
\begin{proof}
  The first part of the theorem follows immediately from the case of
  $sl_2$, see \fullref{sl2 fertile}.
  
  To show the second part we show that roots of $\bs y^{(i)}$ satisfy
  system \eqref{Bethe eqn}, where $\la$ is changed to $s_i\cdot\la$.
  Let $\tilde t_j^{(i)}$ denote the roots of $\tilde y_i$.
  
  The equations of system \eqref{Bethe eqn} corresponding to
  coordinates $\tilde t^{(i)}_j$ are satisfied by \fullref{sl2
    rep}.
  
  The equations of system \eqref{Bethe eqn} corresponding to coordinates
  $\tilde t_j^{(k)}$ such that $a_{ki}=0$, are satisfied because these
  equations are the same for $\tilde{\bs y}$ and $\bs y$.
  
  For any $k$, such that $k \neq i$ and $a_{ik}\neq 0$, choose a root
  $t^{(k)}_j$ of the polynomial $y_k$.  Setting $x=t_j^{(k)}$ in the
  $i$-th equation of \eqref{wr}, we get
\be
\sum_{s}\frac{1}{t_j^{(k)}-t_s^{(i)}}= \sum_{s}\frac{1}{t_j^{(k)}-\tilde
  t_s^{(i)}}+\frac{\langle \la+\rho,\al_i^\vee \rangle}{t_j^{(k)}}.
\ee
This implies that the equation of system \eqref{Bethe eqn}
corresponding to the coordinate $\tilde t^{(k)}_j$ is satisfied as
well.
\end{proof}

Note that if the tuple $\bs y^{(i)}$ is off-diagonal, then it is again
fertile and we can find the $r$--tuple of polynomials $\bs y^{(i,j)}$.
However, in general, we do not know if the tuple $\bs y^{(i)}$ is
off-diagonal.  It is true in the case of $sl_2$ and almost all
non-integral weights $\la$, see \fullref{sl2 rep}.  We have the
following conjecture.

\begin{conj}\label{fer conj}
  An off-diagonal fertile tuple is super-fertile.
\end{conj}

We prove this conjecture for the simple Lie algebras of types $A_r,B_r$,
see Theorems \ref{super} and \ref{super B}.

For an $r$--tuple of polynomials $\bs y$ and a $\g$ weight $\la$, let
$P(\bs y,\la)$ be the set of all pairs of the form $(\bs
y^{(i_1,i_2,\dots,i_m)}, (s_{i_m} \dots s_{i_2}s_{i_1})\cdot \la)$,
where $i_k\in\{1,\dots,r\}$, $m\in\Z_{\geq 0}$, $k=1,\dots,m$, and
$\bs y^{(i_1,i_2,...,i_m)}$ is a descendent of $\bs y$ with respect to
$\la$ in directions $(i_1,\dots,i_m)$.

We call the set $P(\bs y,\la)$ {\it the prepopulation originated at
  $(\bs y,\la)$}.  If $\bs y$ is a super-fertile $r$--tuple with
respect to $\la$ then we call the set $P(\bs y,\la)$ {\it the
  population originated at $(\bs y,\la)$}.

\begin{lem}
  Let $\bs y$ be super-fertile with respect to $\la$ and let $P$ be
  the population originated at $(\bs y,\la)$. Let $(\tilde {\bs
    y},\tilde \la)\in P$. Then $\tilde {\bs y}$ is super-fertile with
  respect to $\tilde\la$ and the population originated at $(\tilde{\bs
    y},\tilde \la)$ is also $P$.  In particular, different populations
  do not intersect.
\end{lem}
\begin{proof}
  By \fullref{ss} we obtain that if $\bs y^{(i)}$ is an immediate
  descendent of $\bs y$ in the direction $i$, then $\bs y$ is also an
  immediate descendent of $\bs y^{(i)}$ in the direction $i$. The
  lemma follows.
\end{proof}

We call a weight $\la$ {\it strongly non-integral} if for any element
of the Weyl group $s\in\mc W$ and any $i\in\{1,\dots,r\}$ the number
$\langle s\cdot\la,\al_i^\vee\rangle$ is not an integer.

Note that if $\la$ is strongly non-integral, then,
in particular, the weights $s\cdot\la$, $s\in \mc W$, 
do not belong to the reflection hyperplanes.
Therefore the map $\mc W\to \mc
W\cdot\la$ mapping $w\in\mc W$ to $w\cdot \la$ is bijective.

\begin{lem}\label{unique}
  Let $\la$ be strongly non-integral and let $r$--tuples $\bs y_1$,
  $\bs y_2$ be descendents of an $r$--tuple $\bs y$ with respect to
  $\la$ in the directions $(i_1,\dots,i_m)$. Then the $r$--tuples $\bs
  y_1$, $\bs y_2$ coincide.
\end{lem}
\begin{proof}
  Lemma follows from the corresponding $sl_2$ statement, see \fullref{sl2 rep}.
\end{proof}

Let $\la$ be strongly non-integral. Let $\bs y$ be super-fertile with
respect to $\la$.  Let $P$ be the population originated at $(\bs
y,\la)$.  For $i\in\{1,\dots,r\}$, let $a_i:\ P\to P$ be the map of
the simple reproduction in the $i$-th direction which maps $(\tilde
{\bs y},\tilde \la)$ to $(\tilde {\bs y}^{(i)},s_i\cdot \tilde \la)$.
According to \fullref{unique}, the map $a_i$ is well defined. By
Lemma $\ref{ss}$ we have $a_i^2=id$. In particular $a_i$ are
invertible.
 
Let $\mc A$ be the subgroup of the group of all permutations of the
elements in $P$ generated by $a_i$, $i=1,\dots,r$.

\begin{conj}\label{conj} 
  There is an isomorphism of groups $\mc A\to \mc W$ which maps $a_i$
  to $s_i$.
\end{conj}

If $(\tilde {\bs y},\tilde \la)\in P$ then $\tilde {\bs y}$ is a
descendent of $\bs y$ in some directions $(i_1,\dots,i_m)$ and we have
$\tilde \la=w\cdot \la$ for some $w\in\mc W$. Since $\la$ is strongly
non-integral, such $w$ is unique and we have $w=s_{i_m}\dots s_{i_1}$.
This defines a map
\bean\label{tau}
\tau:\qquad \ P \;\; \;\; &\to& \mc W,\notag \\
(\tilde{\bs y}, w\cdot \la)&\mapsto& w.
\eean
Since $\bs y$ is super-fertile, $\tau$ is a surjective map.
\fullref{conj} is true if and only if the map $\tau$ is a
bijection for all populations $P$.

Note that the map $\tau=\tau(\bs y, \la)$ depends on the choice of the
element $(\bs y,\la)\in P$. However, if it is bijective for one
element of the population, then it is clearly bijective for all
elements of this population.

\fullref{conj} for the case of $sl_2$ is proved in \fullref{ss}.

Below we prove \fullref{conj} for simple Lie algebras.

\section[Proof of \ref{conj} for simple Lie
  algebras]{Proof of \fullref{conj} for simple Lie
  algebras}\label{proof sec}

\subsection{The case of $sl_{N+1}$}
We have roots $\al_1,\dots,\al_{N}$ with scalar products
$(\al_i,\al_i)=2$, $(\al_i,\al_{i\pm 1})=-1$ and $0$ otherwise.

We fix weights $\bs \La=(\La_1,\dots,\La_n),\La_\infty$, points $\bs
z=(z_1,\dots,z_n)$. The weights $\La_s$, $s=1,\dots,n$, are dominant
integral $sl_{N+1}$ weights and the points $z_s$, $s=1,\dots,n$, are
non-zero, pairwise different complex numbers.  We define polynomials
$T_i$ as in \eqref{T}.

We also fix a strongly non-integral $sl_{N+1}$ weight $\la$.

For any $N$--tuple of functions $\bs y=(y_1,\dots, y_N)$ and $sl_{N+1}$
weight $\la$, we set $y_{N+1}=1$ and define the linear differential
operator of order $N+1$:
\bean\label{Diff}
D (\bs y,\la)
& =& 
(   \partial  - 
\ln'  ( \frac { \prod_{s=1}^N x^{( \la,\al_s)} T_s } { y_N }  ) )\
\dots ( \partial  -  \ln'  ( \frac { y_2 T_1x^{(\la,\al_1)} }{ y_1 }  ) ) \
( \partial - \ln' (  y_1 ) )\notag \\
& =& 
 \prod^{0\ \to\ N}_i\left(
\partial \ -\ \ln' \ \left( 
\frac{ y_{N+1-i}  \prod_{s=1}^{N-i} x^{( \la,\al_s)}
  T_s }{ y_{N-i} }\ \right) \right).
\eean
For an $N$--tuple of polynomials $\bs y$ and an $sl_{N+1}$ weight
$\la$, let as before the prepopulation $P=P(\bs y,\la)$ be the set of
all descendents of $\bs y$ paired with the corresponding weight.

\begin{lem}\label{exist}
  Let $\bs y$ represent an $sl_{N+1}$ critical point associated to
  $(\La_\infty,\la)$. Then the prepopulation $P$ contains the elements
\be
({\bs y}^{(i,i-1,\dots,1)},(s_i s_{i-1} ...s_1)\cdot\la),
\ee 
where $i=0,\dots,N$.
\end{lem}
\begin{proof}
  Since $\bs y$ represents a critical point, it is fertile and none of
  $y_j$ has multiple roots. Moreover we have the Bethe Ansatz equation
  \eqref{Bethe eqn} for each root of each polynomial $y_j$.

In particular there exist polynomials $\tilde y_i$ such that
\be
W(y_i,x^{(\la+\rho,\al_i)}\tilde y_i)=x^{(\la,\al_i)}y_{i-1}y_{i+1}T_i.
\ee
Note that if $t_j^{(i-1)}$ is a root of $y_{i-1}$ then either $\tilde
y_i(t_j^{(i-1)} )\neq 0$ or $\tilde y_i$ vanishes at $t_j^{(i-1)}$ to
order $2$. In the former case, in the same way as in \fullref{gen
  rep}, we see that the Bethe Ansatz equation for the root
$t_j^{(i-1)}$ of $y_{i-1}$ in the $N$--tuple $\bs y^{(i)}$ is still
valid.  In addition the Bethe Ansatz equations for roots of
$y_1,\dots,y_{i-2}$ in the $N$--tuple $\bs y^{(i)}$ are also satisfied
since they are exactly the same as in the $N$--tuple $\bs y$.

Consider the next equation for $\tilde y_{i-1}$: 
\be
W(y_{i-1},x^{(s_i\cdot \la+\rho,\al_{i-1})}\tilde
y_{i-1})=x^{(s_i\cdot\la,\al_{i-1})} y_{i-2}\tilde y_iT_{i-1}.  
\ee 
We have
\be x^{(s_i\cdot \la+\rho,\al_{i-1})}\tilde y_{i-1}=y_{i-1}\int
\frac{x^{(s_i\cdot\la,\al_{i-1})} y_{i-2}\tilde y_iT_{i-1}}{y_{i-1}^2}dx.
\ee 
We claim that the integrand does not have residues. Indeed, the
residues could occur only at the roots $t_j^{(i-1)}$ of $y_{i-1}$. If
$\tilde y_i(t_j^{(i-1)})=0$ then the integrand is holomorphic at
$x=t_j^{(i-1)}$. If $\tilde y_i(t_j^{(i-1)})\neq 0$ then the absence
of the residue is equivalent to the Bethe Ansatz equation
corresponding to $t_j^{(i-1)}$ which is satisfied.  Therefore $\tilde
y_{i-1}$ is a polynomial and there exists a descendent $y^{(i,i-1)}$.

Similarly, we prove that $\bs y^{(i,i-1,\dots,i-m)}$ is a well-defined
$N$--tuple of polynomials for $m=2,\dots,i-1$.
\end{proof}

If $(\tilde {\bs y},\tilde \la)\in P$ then we write
\be
\tilde \la=w(\tilde \la)\cdot \la=\la-\sum_{i=1}^N a_i(\tilde \la)\al_i,
\ee
where $w(\tilde \la)\in\mc W$ and $a_i(\tilde \la)\in \R$.  Define
{\it the shifted prepopulation} $\bar P(\bs y,\la)$ as the following
set of $N$--tuples of functions: 
\be \bar P=\{(x^{a_1(\tilde
  \la)}\tilde y_1,x^{a_2(\tilde \la)}\tilde y_2,\dots, x^{a_N(\tilde
  \la)}\tilde y_N) \ | \ (\tilde{\bs y},\tilde \la) \in P\}.  
\ee 
Note that the shifted prepopulation $\bar P(\bs y,\la)$ depends on the
choice of an element $(\bs y,\la) \in P$. However, the difference is
not very essential: if $(\tilde y,\tilde \la)\in P$ and $\bar P(\tilde
y,\tilde \la)$ is the corresponding shifted prepopulation then there
exists $a_i\in\R$ such that the $i$-th function in any $N$--tuples in
$\bar P(\tilde y,\tilde \la)$ is obtained via multiplication of the
$i$-th function of the corresponding $N$--tuple in $\bar P(\bs y,\la)$
by $x^{a_i}$.

\begin{lem}\label{ind}
  Let $\tilde{\bs y}\in\bar P(\bs y,\la)$. 
  Then $D(\tilde {\bs  y},\la)=D(\bs y,\la)$.
\end{lem}
\begin{proof}
  Let $(\bs v,\mu)\in P$ and let $(v^{(i)},s_i\cdot \mu)\in P$ be the
  immediate descendent of $\bs v$ with respect to $\mu$ in the
  direction $i$.  Then we have $v^{(i)}_k=v_k$ for $k\neq i$ and
  $W(v_i,v^{(i)}_ix^{(\mu+\rho,\al_i)})=x^{(\mu,\al_i)}
  v_{i-1}v_{i+1}T_i$. The last relation can be rewritten as
\bean\label{try}
W(x^{a_i(\mu)}v_i,v^{(i)}_ix^{a_i(s_i\cdot \mu)})=
x^{(\la,\al_i)}\ v_{i-1}x^{a_{i-1}(\mu)} \ v_{i+1}
x^{a_{i+1}(\mu)} T_i.
\eean
Let 
\bea
(\bar{\bs v},\mu) &=
& ((x^{a_1(\mu)}v_1,\dots,x^{a_N(\mu)}v_N),\mu)\in \bar P(\bs y,\la),\\
(\bar{\bs v}^{(i)},s_i\cdot \mu) &=&
((x^{a_1(s_i\cdot\mu)}v_1^{(i)},\dots,
x^{a_N(s_i\cdot\mu)}v_N^{(i)}),s_i\cdot\mu)\in
\bar P(\bs y,\la).
\eea
Identity \eqref{try} reads $W(\bar v_i,\bar v_i^{(i)})=x^{
 ( \la,\al_i)}\bar v_{i-1} \bar v_{i+1}T_i$ and therefore
\bean\label{same middle}
-\ln'' (\bar v_i)+\ln'(\bar v_i)(\ln'( T_i x^{(\la,\al_i)}
\bar v_{i-1} \bar v_{i+1}) -\ln'(\bar v_i))=\notag \\
-\ln'' (\bar v_i^{(i)})+\ln'(\bar v_i^{(i)})(\ln'( T_i
x^{(\la,\al_i)} 
\bar v_{i-1}^{(i)}\bar v_{i+1}^{(i)}) -\ln'(\bar v_i^{(i)})).
\eean
Compare $D(\bar {\bs v})$ with $D(\bar {\bs v}^{(i)})$.  All factors
but two successive ones in these operators are the same. The products
of the two middle factors are the same by \eqref{same middle}
\end{proof}

We call $D(\bs y,\la)=D$ {\it the operator associated to the shifted
  prepopulation $\bar P(\bs y,\la)$}.  It follows from \fullref{ind}
if $\tilde {\bs y}\in\bar P(\bs y,\la)$ then $D\tilde y_1=0$.

We use the following notation for Wronskians and {\it divided Wronskians}:
\bea
W(u_1,\dots,u_i)& =& \det(u_k^{(j-1)})_{k,j=1}^i,\\
W^\dagger(u_1,\dots,u_i) &=& \frac{W(u_1,\dots,u_i)}{(x^{(
    \la,\al_1)}T_1)^{i-1}(x^{(
    \la,\al_2)}T_2)^{i-2}\dots (x^{(
    \la,\al_{i-1})} T_{i-1})},
\eea
$i=1,\dots,N+1$, where $u_k^{(j-1)}$ denotes the $(j-1)st$ derivative
of $u_k$ with respect to variable $x$.

\begin{lem}\label{basis} 
  Let $\bs y$ either represent an $sl_{N+1}$ critical point associated
  to $(\La_\infty,\la)$ or be super-fertile with respect to $\la$.
  Then there exist functions $u_1,\dots, u_{N+1}$ such that $Du_i=0$,
  $y_i=W^\dagger(u_1,\dots,u_i)$ for $i=1,\dots,N+1$ and
  $u_ix^{-(\la+\rho,\al_1+\dots +\al_{i-1})}$ are polynomials for
  $i=1,\dots,N+1$.
\end{lem}
\begin{proof}
  Let $u_i$ to be the first coordinate of the $N$--tuple in the element
  of $\bar P(\bs y,\la)$ corresponding to ${\bs
    y}^{(i-1,i-2,\dots,1)}$ which is a descendent of $\bs y$ with
  respect to $\la$ in directions $(i-1,i-2,\dots,1)$.  If $\bs y$
  represents an $sl_{N+1}$ critical point with respect to $\la$ then
  such an $N$--tuple exists by \fullref{exist}.
  
  We have $Du_i=0$ by \fullref{ind}.
  
  The condition $W^\dagger(u_1,\dots,u_i)= y_i$ follows from the
  standard Wronskian identities, cf \cite[proof of Lemma 5.5]{MV}.
\end{proof}

\begin{cor}\label{sing points}
  Let $\bs y$ be an $N$--tuple of polynomials and $l_i=\deg y_i$, 
  $i=1,\dots,N$. Let $\bs \La_\infty$ be given by \eqref{at inf}.
  Let $\bs y$ represent an $sl_{N+1}$ critical point associated to
  $(\La_\infty,\la)$ or let $\bs y$ be super-fertile with respect to
  $\la$.  Then the kernel of
  the operator $D(\bs y,\la)$ is spanned by functions of the form
\bean\label{kernel form}
p_0, p_1x^{(\la+\rho,\al_1)}, \dots, p_Nx^{(\la+\rho,\al_1+\dots+\al_N)},
\eean
where $p_i$ is a polynomial of degree $\deg
y_1+(\La_\infty,\al_1+\dots+\al_i)$, $p_i(0)\neq 0$.  The only
singular points of the operator $D(\bs y,\la)$ in $\C^*$ are regular
singular points located at $z_1,\dots,z_n$, and the exponents at
$z_i$, $i=1,\dots,n$, are 
\bean\label{exponents} 
z_i: 0,(\La_i+\rho,\al_1),
(\La_i+\rho,\al_1+\al_2),\dots,(\La_i+\rho,\al_1+\dots+\al_N), 
\eean
Conversely, if a linear differential operator $D$ of order $N+1$ has
the kernel spanned by functions of the form \eqref{kernel form} and the
only non-zero singular points of $D$ in $\C^*$ are regular singular
points at $z_1,\dots,z_n$ with the exponents given by \eqref{exponents},
then the $N$--tuple $\bs y$ given by the divided Wronskians
\be
y_i=
W^\dagger(p_0, p_1x^{(\la+\rho,\al_1)}, \dots, 
p_{i-1}x^{(\la+\rho,\sum_{j=1}^{i-1}\al_j)})
\ee
is super-fertile with respect to $\la$ and satisfies 
$\deg y_i=l_i$, $i=1,\dots,N$.
\end{cor}
\begin{proof}
 For $i=0,\dots,N$, we set 
\be
p_i=u_{i+1}x^{-(\la+\rho,\al_1+\dots +\al_{i})},
\ee
where $u_1,\dots,u_{N+1}$ are as in \fullref{basis}. 
By \fullref{basis}, $p_0,\dots,p_N$ are polynomials.
 
 Now, the first part of \fullref{sing points} 
 follows from \fullref{basis} by standard Wronskian identities, cf
 \cite[Lemmas 5.8 and 5.10]{MV}.

Conversely, let $V$ be the kernel of the operator $D$.
  We have $(N+1)!$ distinguished full flags in $V$ such that the
  divided Wronskians of all spaces which form the flags are of the
  form $x^ap(x)$ where $p(x)$ is a polynomial. Namely, for a
  permutation $w\in S_{N+1}$ of the set $\{0,1,\dots,N\}$ 
 we have a full flag $\mc F_w$ such that
  the space of dimension $i$ is spanned by 
  $p_{w(0)}x^{(\la+\rho),\al_{w(0)}}\dots,
  p_{w(i-1)}x^{(\la+\rho),\al_{w(i-1)}}$.
  
  For each such flag $\mc F_w$, we have the corresponding element $p_w$ in
  $\bar P(\bs y,\la)$. Each element $p_w$ is obviously fertile and the
  immediate descendents of $p_w$ in the $i$-th direction is
  $p_{(i,i+1)w}$. Therefore, the $N$--tuple $\bs y$ is fertile with 
   respect to $\la$.
\end{proof}

\begin{theorem}\label{super}
  Conjectures \ref{fer conj}, \ref{conj} hold for the case of
  $sl_{N+1}$.
\end{theorem}
\begin{proof}
  \fullref{fer conj} follows from \fullref{gen rep} 
and \fullref{sing points}.
  
\fullref{conj} in the case of $sl_{N+1}$
follows from the proof of the converse statement of  
\fullref{sing points}.
\end{proof}

\subsection{The case of $B_2$}
In the case of $B_2$ we have two roots $\al_1,\al_2$ such that
$(\al_1,\al_1)=4$, $(\al_2,\al_2)=2$, $(\al_1,\al_2)=-2$.

The key observation is that $B_2$ populations can be embedded in
$sl_4$ populations.

Given a $B_2$ weight $\La$, define the $sl_4$ weight $\La^A$ by
\be
\langle \La^A , (\al_1^A)^\vee \rangle = \langle \La^A , (\al_{3}^A)^\vee
\rangle  =
\langle \La  ,  \al_1^\vee \rangle , \qquad \langle \La^A , (\al_{2}^A)^\vee
\rangle  =
\langle \La  ,  \al_2^\vee \rangle,
\ee
where $\al_i^A$ are roots of $sl_4$.

Note that if $\La$ is strongly non-integral then $\La^A$ is strongly
non-integral.

\begin{lem} A pair $\bs y=(y_1,y_2)$ represents a $B_2$ critical point
  associated to $(\bs z,\bs \La$; $\La_\infty,\la)$ if and only if the
  triple $\bs y^A=(y_1,y_2,y_1)$ represents an $sl_4$ critical point
  associated to $(\bs z,\bs \La^A$; $\La_\infty^A,\la^A)$. Moreover
  there is an embedding $P(\bs y,\la) \to P^A(\bs y^A,\la^A)$ which
  sends $((\tilde y_1,\tilde y_2),\tilde \la)\in P(\bs y,\la)$ to
  $((\tilde y_1,\tilde y_2,\tilde y_1),\tilde \la^A)\in P^A(\bs
  y^A,\la^A)$.
\end{lem}
\begin{proof}
  Follows immediately from the definitions.
\end{proof}

\begin{theorem}\label{b2}
  \fullref{conj} holds in the case of root system $B_2$.
\end{theorem}
\begin{proof}
  Recall the surjective maps $\tau: P\to \mc W$ and $\tau^A: P^A\to \mc
  W^A$, where $\mc W$ and $\mc W^A$ are the $B_2$ and $A_3$ Weyl
  groups, see \eqref{tau}.  Then we clearly have $(\tau(\bs
  y)\cdot\la)^A=\tau^A(\bs y^A)\cdot \la^A$.  By \fullref{super},
  the map $\tau^A$ is injective.
  
  Therefore $\tau$ is injective and hence bijective.
\end{proof}

\begin{theorem}\label{super B}
  \fullref{fer conj} holds in the case of simple Lie algebras of type
  $B_N$.
\end{theorem}
\begin{proof}
  Similarly to the case $N=2$, the $N$--tuple $(y_1,\dots,y_N)$
  represents a critical point of type $B_N$ if and only if the
  $(2N-1)$--tuple $(y_1,\dots,y_{N-1},y_N,y_{N-1},\dots,y_1)$
  represents a critical point of type $A_{2N-1}$, see also \cite{MV}.
  The tuple $(y_1,\dots,y_N)$ is super-fertile in $B_N$ sense because
  the $(2N-1)$--tuple $(y_1,\dots,y_{N-1},y_N,y_{N-1},\dots,y_1)$ is
  super-fertile in $A_{2N-1}$ sense.
\end{proof}

\subsection{The case of the root systems of types $B,C,D,E$ and $F$}
From the $sl_2$, $sl_3$ and $B_2$ cases, we obtain the general case
(except for $G_2$).

\begin{theorem}\label{general}
  \fullref{conj} holds in the case of the root systems of types
  $B,C,D,E$ and $F$.
\end{theorem}
\begin{proof}
  Let $\g$ be a rank $r$ simple Lie algebra of type $B,C,D,E$ or $F$.
  
  The Weyl group of $\g$ is a finite group described by the generators
  and relations: 
\be 
\mc W=<s_1,\dots,s_r>/(s_i^2=(s_is_j)^{-(\al_i,\al_j)+2}=1,\ 
  i,j=1,\dots,r,\ i\neq j).  
\ee 
  Here $<s_1,\dots,s_r>$ denotes the
  free group with generators $s_1,\dots,s_r$.
  
  Note that in our case $(\al_i,\al_j)$ takes values $2,0,-1$ or $-2$.
  
  Let $\la$ be strongly non-integral, $\bs y$ represent a
  $\g$--critical point associated to $(\La_\infty,\la)$ and let
  $P=P(\bs y,\la)$ be the corresponding population.
  
  We have the corresponding relations in the populations among
  descendents:
\begin{itemize}
\item  $\bs y^{(i,i)}=\bs y$, 
\item  $\bs y^{(i,j)}=\bs y^{(j,i)}\ $ if
  $(\al_i,\al_j)=0$,
\item  $\bs y^{(i,j,i)}=\bs y^{(j,i,j)}\ $ if
  $(\al_i,\al_j)=-1$, 
\item $\bs y^{(i,j,i,j)}=\bs
  y^{(j,i,j,i)}\ $ if
  $(\al_i,\al_j)=-2$. 
\end{itemize}
The first relation follows from the case of $sl_2$, \fullref{ss},
the second relation is obvious, the third relation follows from the
case of $sl_3$, \fullref{super}, and the fourth relation follows
from the case of $B_2$, \fullref{b2}.

Therefore, $P$ has at most $|\mc W|$ elements.
\end{proof}

\subsection{The case of $G_2$}
In the case of $G_2$ we have two roots $\al_1,\al_2$ such that
$(\al_1,\al_1)=2$, $(\al_2,\al_2)=6$, $(\al_1,\al_2)=-3$.

The key observation is that $G_2$ populations can be embedded in $C_3$
populations.

Given a $G_2$ weight $\La$, define the $C_3$ weight $\La^C$ by
\be
\langle \La^C , (\al_1^C)^\vee \rangle = \langle \La^C , (\al_{3}^C)^\vee
\rangle  =
\langle \La  ,  \al_1^\vee \rangle , \qquad \langle \La^C , (\al_{2}^C)^\vee
\rangle  =
\langle \La  ,  \al_2^\vee \rangle,
\ee
where $\al_i^C$ are roots of $C_3$.

Note that if $\La$ is strongly non-integral then $\La^C$ is strongly
non-integral.

\begin{lem} A pair $\bs y=(y_1,y_2)$ represents a $G_2$ critical point
  associated to $(\bs\La,\bs z$; $\La_\infty,\la)$ if and only if the
  triple $\bs y^C=(y_1,y_2,y_1)$ represents a $C_3$ critical point
  associated to $(\bs \La_i^C,\bs z$; $\La_\infty^C,\la^C)$.  Moreover
  there is an embedding $P(\bs y,\la) \to P^C(\bs y^C,\la^C)$ which
  sends $((\tilde y_1,\tilde y_2),\tilde \la)\in P(\bs y,\la)$ to
  $((\tilde y_1,\tilde y_2,\tilde y_1),\tilde \la^C)\in P^C(\bs
  y^C,\la^C)$.
\end{lem}
\begin{proof}
  Follows immediately from definitions.
\end{proof}

\begin{theorem}\label{g2}
  \fullref{conj} holds in the case of the root system $G_2$.
\end{theorem}
\begin{proof}
  Recall the surjective maps $\tau: P\to \mc W$ and $\tau^C: P^C\to
  \mc W^C$, where $\mc W$ and $\mc W^C$ are the $G_2$ and $C_3$ Weyl
  groups, see \eqref{tau}.  Then we clearly have $(\tau(\bs
  y)\cdot\la)^C=\tau^C(\bs y^C)\cdot \la^C$.  By \fullref{general}, the map $\tau^C$ is injective.
\end{proof}

To summarize, we have
\begin{cor}\label{general cor}
  \fullref{conj} holds in the case of all simple Lie algebras.
\end{cor}

\section{Weyl group actions on Bethe vectors}\label{act sec}

\subsection{The trigonometric Gaudin operators and Bethe vectors}
Let $\g$ be a simple complex Lie algebra of rank $r$ with the Killing 
form $(\ ,\ )$. 
We choose a Cartan subalgebra $\h$, simple roots $\al_i$, $i=1,\dots,r$. 
We identify $\h$ with $\h^*$ using the Killling form on $\g$. 
Let $F_i$, $E_i$, $i=1,\dots, r$, be the Chevalley generators of $\g$. 

Let $\Delta$ be the root system, 
let $\Delta_\pm$ be the sets of positive and negative roots and 
let $\g=\h\oplus\sum_{\al\in\Delta}g_\al$ be the root
decomposition. Let $e_\al$, $\al\in\Delta$, 
be generators of $g_\al$ such that $(e_\al,e_{-\al})=1$. 
Let $\{h_j\}_{j=1,\dots, r}$ be an orthonormal basis of the Cartan
algebra $\h\subset \g$. 

Set 
\be
\Om^0=\frac12\sum_{j=1}^r h_j\otimes h_j, \qquad \Om^+=\Om^0+\sum_{\al\in\Delta_+}e_\al\otimes e_{-\al},\qquad \Om^-=
\Om^0+\sum_{\al\in\Delta_+}e_{-\al}\otimes e_{\al}.
\ee
The {\it trigonometric $R$--matrix} is defined by
\be
r(z)=\frac{\Om^+z+\Om^-}{z-1}.
\ee


We fix $\bs z, \bs \La,\La_\infty,\bs l,\la$ as in \fullref{master
  sec}. Let $L_1,\dots, L_n$ be irreducible $\g$--modules with highest weights
$\La_1,\dots, \La_n$ and let $V=L_1\otimes \dots \otimes L_n$.
Let $V[\mu]\subset V$ be the subspace of $V$ of all vectors of weight $\mu$.

We write $X^{(k)}$ for an operator $X\in\g$
acting on the $k$-th factor. Similarly we write $X^{(k,l)}$ for an
operator $X\in\g\otimes\g$ acting on the $k$-th and $l$-th factors.

The {\it trigonometric Gaudin
operators} $H_i(\la)$, $i=1,\dots,n$, are defined by
\be
H_i(\la)=\la^{(i)}+\sum_{j=1,\dots,n, \ j\neq i} r^{(i,j)}(z_i/z_j).
\ee
The trigonometric Gaudin operators depend on $\la\in\h$ and act in $V$.
The trigonometric Gaudin operators all commute, $[H_i(\la),H_j(\la)]=0$, $i,j=1,\dots,n$.
The trigonometric Gaudin operators commute with the action of $\h$ on $V$ and 
therefore preserve every weight subspace of $V$.

For a given $\la$, common eigenvectors of the trigonometric 
Gaudin operators $H_i(\la)$ in the weight subspace $V[\Lambda_\infty]$
can be constructed by the Bethe Ansatz method as follows.

Let $l=l_1+\dots+l_n$.  Let $c$ be the unique non-decreasing function
from $\{1,\dots,l\}$ to $\{1,\dots,r\}$ such that $\sharp\
c^{-1}(i)=l_i$, $i=1,\dots,r$.

Let $P(\bs l,n)$ be the set of sequences $I\ = \ (i_1^1, \dots ,
i^1_{j_1};\ \dots ;\ i^n_1, \dots , i^n_{j_n})$ of integers in $\{1,
\dots , r\}$ such that for all
$i = 1, \dots ,  r$, the integer $i$ appears in $I$ precisely $l_i$ times. 
For $I \in P(\bs l, n)$, and a permutation $\sigma \in S_l$
set $\sigma_1(i) = \sigma(i)$ for $i = 1, \dots , j_1$,
and $\sigma_s(i) = \sigma(j_1+\cdots +j_{s-1}+i)$ for $s = 2, \dots , n$ and 
$ i = 1, \dots , j_s$. Define 
\bea
S(I)\ {} = \ {}
\{\ \sigma \in S_l\ {} |\ {} c(\sigma_s(j)) = i_s^j 
\ {} \text{for} \ {} s = 1, \dots , n \ {} \text{and} \ {} 
j = 1, \dots , j_s\ \}\ .
\eea 

For $I \in P(\bs l, n)$ we define a vector in $V[\La_\infty]$ by the formula
\bea
F_Iv\ =\ F_{i_1^1} \dots F_{i_{j_1}^1}v_{1} \otimes \cdots
\otimes F_{i_1^n} \dots F_{i_{j_n}^n}v_{n}.
\eea

For $I \in P(\bs l, n)$, $\sigma\in S(I)$, we define 
a rational function of $\bs t=(t_i^{(j)})_{j=1,\dots,r}^{i=1,\dots,l_j}$ 
by the formula
\bea
\omega_{I,\sigma}(\bs t) \ =\ 
\omega_{\sigma_1(1),\ldots,\sigma_1(j_1)}(z_1;\bs t)\
\cdots\
\omega_{\sigma_n(1),\ldots,\sigma_n(j_n)}(z_n;\bs t) ,
\eea
where 
\bea
\omega_{i_1,\ldots, i_j}(z;\bs t) \ =\ \frac 1 {(t_{i_1}-t_{i_2}) \cdots
(t_{i_{j-1}}-t_{i_j}) (t_{i_j}-z)}\ 
\eea
and
$(t_1,\dots,t_l)=
(t_1^{(1)},\dots,t_{l_1}^{(1)},t_1^{(2)},
\dots,t_{l_2}^{(2)},\dots,t_1^{(r)},\dots,t_{l_r}^{(r)})$.

We define {\it the weight function} by 
\bean\label{bethe vector}
\omega(\bs t) \ =\  \sum_{I\in P(\bs l,n)}\ \sum_{\sigma\in S(I)}\ 
\omega_{I,\sigma}(\bs t)\ F_I v\ .
\eean
The weight function $\om(\bs
t)$ is a rational function of $\bs t$ with values in  $V[\La_\infty]$.

\begin{prop}\label{Bethe vector prop} 
  If $\bs t$ is a critical point of the master function
  \eqref{master} associated to $(\La_\infty, \la)$, then
  $\om(\bs t)$  is a well defined vector of weight $\La_\infty$ in $V$ 
   which is an eigenvector of the operators 
  $H_i(\la+\rho+\La_\infty/2)$, $i = 1,\dots,n$.
\end{prop}
\begin{proof} 
  \fullref{Bethe vector prop} follows from the corresponding
  fact for the rational Gaudin operators, see \cite{SV} and the
  relation between rational and trigonometric Gaudin operators, see
  \cite[Appendix B]{MaV}.
\end{proof}

The vector $\om(\bs t)$ is called the {\it Bethe vector associated to
  the critical point $\bs t$}. It is expected that for generic values
of parameters, all critical points are non-degenerate and the Bethe
vectors form a basis in $V$. In particular, the number of orbits of
critical points and thus the number of populations should match the
dimension of the subspace of all vectors of weight $\La_\infty$ in $V$.

\subsection{Counting $sl_{N+1}$ critical points}
Let $L_\La$ denote the irreducible $sl_{N+1}$ module of highest weight $\La$.
\begin{prop}\label{sln number}
  For almost all $\la$ the number of orbits of $sl_{N+1}$ critical
  points associated to $(\La_\infty,\la)$ and counted with
  multiplicity does not exceed the dimension of the subspace of the
  weight $\La_\infty$ in the tensor product
  $L_{\La_1}\otimes\dots\otimes L_{\La_n}$.
\end{prop}
\begin{proof}
  By \fullref{finite}, for almost all $\la$, the number of critical
  points associated to $(\La_\infty,\la)$ is finite. Therefore, there
  is a Zariski open set $O\subset \bar \h^*$, such that the number of
  orbits of $sl_{N+1}$ critical points associated to
  $(\La_\infty,\la)$ and counted with multiplicities is the same for
  all $\la\in O$.
  
  If $\la$ is a dominant integral weight, then the number of orbits of
  critical points associated to $(\La_\infty,\la)$ and counted with
  multiplicities is bounded from above by the multiplicity of
  $L_{\La_\infty}$ in the tensor product $L_\la\otimes
  L_{\La_1}\otimes\dots\otimes L_{\La_n}$, see \cite{MV,BMV}.
  
  For any integer $M>0$, let $\mc C_M$ be the set of all weights
  $\la\in\h^*$ such that the scalar products $(\la,\al_i)$ are
  integers greater than $M$.
  
  If $\la\in \mc C_M$ and $M$ is large enough, then any singular vector of weight $\La_\infty $ in the tensor product $L_\la\otimes
  L_{\La_1}\otimes\dots\otimes L_{\La_n}$ is uniquely determined by its projection to $v_\la\otimes L_{\La_1}\otimes\dots\otimes L_{\La_n}$, 
 where $v_\la$ is the highest weight vector of $L_\la$.
  Therefore, the multiplicity
  of $L_{\La_\infty}$ in the tensor product $L_\la\otimes
  L_{\La_1}\otimes\dots\otimes L_{\La_n}$ equals the dimension of the
  subspace of weight $\La_\infty$ in the tensor product
  $L_{\La_1}\otimes\dots\otimes L_{\La_n}$.

  For any $M\in\Z_{\geq 0}$, the set $\mc C_M$ is not contained in any
  proper algebraic subset in $\bar \h^*$ and thus the proposition follows.
\end{proof}

\begin{prop}\label{sl2 number}
  For almost all $\la$ and almost all $(z_1,\dots,z_n)\in\C^n$, all of
the critical points are non-degenerate and the number of orbits of
$sl_2$ critical points associated to $(\La_\infty,\la)$ equals the
dimension of the subspace of weight ${\La_\infty}$ in the tensor
product $L_{\La_1}\otimes\dots\otimes L_{\La_n}$.
\end{prop}
\begin{proof}
  For almost all $\la$ the number of orbits of $sl_2$ critical points
  associated to $(\La_\infty,\la)$ is the same.

  If $\la$ is dominant integral, then the number of orbits of critical points
  associated to $(\La_\infty,\la)$ for generic $\bs z$ equals the
  multiplicity of $L_{\La_\infty}$ in the tensor product $L_\la\otimes
  L_{\La_1}\otimes\dots\otimes L_{\La_n}$ and all these orbits are
  non-degenerate, see \cite[Theorem 1]{ScV}.
  
  For dominant integral values of $\la$ which are large enough, the multiplicity
  of $L_{\La_\infty}$ in the tensor product $L_\la\otimes
  L_{\La_1}\otimes\dots\otimes L_{\La_n}$ equals the dimension of the
  subspace of weight $\La_\infty$ in the tensor product
  $L_{\La_1}\otimes\dots\otimes L_{\La_n}$.
  
  In this case, almost all $\la$ means all but finitely many and
  therefore the proposition follows.
\end{proof}

\subsection{Actions of the Weyl group on Bethe vectors}
Let $\g$ be a simple Lie algebra, $G$ the corresponding connected and
simply connected Lie group.  The group $G$ acts on any
finite-dimensional irreducible representation of $\g$.  Let $\h\subset
\g$ be a Cartan subalgebra, and $T\subset G$ the corresponding torus.
The Weyl group of $\g$ can be described as $N/T$ where
$N=\{g\in G\ |\ gTg^{-1}= T\}$. In particular, this defines a
projective action of the Weyl group on any tensor product of
finite-dimensional irreducible representations of $\g$. The projective
action becomes an action in the zero weight subspace.

We fix our $\bs z, \bs \La,\La_\infty,\bs l,\la$ as in \fullref{master
  sec}. Let $L_1,\dots, L_n$ be irreducible $\g$--modules with highest weights
$\La_1,\dots, \La_n$ and let $V=L_1\otimes \dots \otimes L_n$. 
Let $V[\mu]$ be the subspace of all vectors in $V$ of weight $\mu$. Let 
$P(V)=\{\mu,\  V[\mu]\neq 0\}$ be the set of all nontrivial weights in $V$. 

We define the dynamical Weyl group acting on $V$ following \cite{TV}.

Let $M_\mu$ denote the Verma module with highest weight $\mu$, \
 $v_\mu$ a highest weight vector in $M_\mu$.

Let $M_\mu, M_\la$ be Verma modules. Two cases are possible:
\begin{enumerate}
\item [a)] 
Hom$_\g(M_\mu,M_\la)=0$ or 
\item [b)]  Hom$_\g(M_\mu,M_\la)=\C$ and every nontrivial homomorphism $M_\mu\to M_\la$
is an embedding.
\end{enumerate}

Let $M_\la$ be a Verma module with dominant weight $\la\in P^+$.
Then Hom$_\g(M_\mu,M_\la)=\C$ if and only if there is $w\in \W$
such that $\mu= w\cdot \la$.

Let $w=s_{i_k}\ldots s_{i_1}$ be a reduced presentation 
of an element of the Weyl group $\mc W$.
Set $\al^{1}=\al_{i_1}$ and $\al^{j}=(s_{i_1}\ldots s_{i_{j-1}})(\al_{i_j})$
for $j=2,\ldots,k$. Let $n_j=(\la+\rho,(\al^{j})^\vee)$.
For a dominant $\la\in P^+$, the numbers $n_j$ are positive integers.
Define a singular vector $v_{w\cdot\la}^\la\in M_\la$ by
\bean
v_{w\cdot\la}^\la\,=\,
{(E_{-\al_{i_k}})^{n_k}\over n_1!}
\ldots {(E_{-\al_{i_1}})^{n_1}\over n_k!}\,v_\la\,.
\eean
This vector does not depend on the
reduced presentation, see \cite{TV}.

For all $\la\in P^+$, $w\in \mc W$, fix an embedding $M_{w\cdot \la}
\hookrightarrow M_\la$ sending $v_{w\cdot \la}$ to
$v_{w\cdot \la}^ \la$.

We say that  $\la\in P^+$ is generic with respect to $V$ if
\begin{itemize}
\item For any $\nu\in P(V)$ and any $v \in V[\nu]$, 
there exist a unique intertwining operator
$\Phi^v_\la:M_\la\to M_{\la-\nu}\otimes V$ such that $\Phi^v_\la (v_\la)=
v_{\la-\nu}\otimes v + $  terms of lower weight in the first factor.
\item
For any $w,w' \in \mc W,\, w\neq w'$, and any $\nu \in P(V)$,
the vector $w\cdot\la -w'\cdot(\la -\nu)$ does not belong to $P(V)$.
\end{itemize}
If $\la=\sum_i \la_i\om_i$, where $\om_i$ are fundamental weights and $\la_i$ are large enough positive numbers then $\la$ is generic with respect to $V$.

\begin{lem}\label{a}{\rm\cite{TV}}\qua
Let $\la\in P^+$ be generic with respect to $V$. Let $v\in V[\nu]$.
Consider the intertwining operator $\Phi^v_\la:M_\la\to M_{\la-\nu}\otimes V$.
For $w\in \W$, consider the singular vector $v_{w\cdot \la}^\la\in M_\la$.
Then there exists a unique vector $A_{w}(\la)(v)\in V[w(\nu)]$ such that
\bean\label{A}
\Phi_\la^v (v_{w\cdot\la}^\la)= v_{w\cdot(\la-\nu)}^{ \la-\nu}\otimes
A_{w}(\la)(v)\,
+\,\text{terms of lower weight in the first factor}\,.
\notag
\eean
\end{lem}

For  generic $\la\in P^+$,
\fullref{a} defines
a linear operator $A_{w}(\la): V\to V$ such that $A_{w}(\la)(V[\nu])
\subset V[w(\nu)]$ for all $\nu\in P(V)$. 
This operator is extended to other values of $\la$ as
a rational function of $\la$.

The collection of rational functions
$A_{w}(\la)$, $w\in\mc W$, is called 
{\it the dynamical Weyl group acting on $V$}.

Introduce new linear 
operators $\mc A_{w}(\la)\, :\, V\,\to\, V$ for $w\in\mc W$. 
Namely, for any
$w\in \mc W,\, \nu \in P(V),\, v\in V[\nu]$, set
$$
\mc A_w (\la) \,v\ = A_w(\la + \nu)\, v \ .
$$
We still have $\mc A_{w}(\la)(V[\nu])
\subset V[w(\nu)]$ for all $\nu\in P(V)$.

\begin{lem}{\rm\cite{TV}}
\label{proper}

\begin{itemize}
\item
For any 
$w_1,w_2 \in \mc W$ and 
$\nu \in P(V)$, we have
\be
\left(
\mc A_{w_1}(w_2 \cdot \la)\,
\mc A_{w_2}(\la)
\right)\,
\vert_{V[\nu]} =
{c}_{w_1,w_2,\la,\nu}
\ \mc A_{w_1w_2}(\lambda)\, \vert_{V[\nu]}\ ,
\ee
where ${c}_{w_1,w_2, \la,\nu}$ is a constant depending on $w_1, w_2,\la, \nu$.

\item
For any $w, w_1,w_2\in \mc W$, $\nu \in P(V)$, the limits 
$$
\mc A_w (\infty) = \lim_{\la\to \infty}\ \mc A_{w}(\la),
\qquad c_{w_1,w_2,\nu}=  \lim_{\la\to \infty}\ c_{w_1,w_2,\la,\nu}
$$ 
do exist. Therefore, we have
\be
(\mc A_{w_1}(\infty)\,
\mc A_{w_2}(\infty))\,\vert_{V[\nu]}\ =\
{c}_{w_1,w_2,\nu}
\ \mc A_{w_1w_2}(\infty)\, \vert_{V[\nu]}.
\ee
Moreover, the collection of operators $\mc A_w(\infty), \, w\in \mc W$, 
gives the canonical projective action of $\mc W$ on $V$.

\item
For any vector $v\in V[\nu]$ and $w\in\mc W$, we have
\be
\mc A_{w}(\la)\ H_i\left(\la+\rho+\frac\nu 2\right)\,v\ =\
H_i\left(w\cdot \la+\rho+\frac{w(\nu)}2\right)\ \mc A_{w}(\la)\,v\ .
\ee
\end{itemize}
\end{lem}
\begin{proof}
The first statement follows from \cite[Theorems 8 and 10]{TV}.
The second statement is \cite[Corollary 14]{TV}.
The statement of \cite[Lemma 18]{TV}, which holds for any root
system, gives the last statement of our lemma.
\end{proof}

According to this lemma, if $\omega$ is an eigenvector of the operators
$H_i(\la+\rho+\frac\nu2)$,   then $\mc A_w(\la)\, \omega$ is an
eigenvector of the operators
$H_i(w\cdot \la+\rho+\frac{w(\nu)}2)$.  

Let $\la$ be generic.  Let $\bs t$ be a solution of the Bethe Ansatz
equation associated to $(\La_\infty,\la)$ and let $\bs y$ be the
corresponding $r$--tuple of polynomials. Then $\om(\bs t)$ is an eigenvector  
of the operators $H_i(\la+\rho+\La_\infty/2)$. 

By \fullref{general  cor}, 
for each element $w$ of the Weyl group, we have the descendent
$w\bs y$ of $\bs y$ obtained via the reproduction procedure.  Let
$w\bs y$ represent the tuple $\bs t_w$. Moreover, if $w\bs y$ is
off-diagonal,\ then $\bs t_w$ is a critical point associated to
$(w\La_\infty,w\cdot\la)$ and $\om(\bs t_w)$ is an eigenvector  
of the operators $H_i(w\cdot \la+\rho+w\La_\infty/2)$.

We conjecture that the action of the operator $\mc A_w(\la)$ 
coincides with the action of 
the Weyl group, induced by the reproduction procedure 
(when the latter action is well-defined). 
More precisely, we have

\begin{conj}\label{Weyl acts}  Let $\la$ be generic.
  Let $\bs t$ be a critical point of the master function \eqref{master}
  associated to $(\La_\infty,\la)$ and let $\bs y$ be the
corresponding $r$--tuple of polynomials. Let $\om(\bs t)\in V[\La_\infty]$ be
  the corresponding Bethe vector.  Let $w \in \mc W$. 
 Assume that $w\bs y$ is off-diagonal. Let $w\bs y$ represent the tuple 
$\bs t_w$.

Then the vector
  $\mc A_{w}(\la)\,\om(\bs t)$ is a scalar multiple of the Bethe
  vector $\om(\bs t_w)$.
\end{conj}
Below we prove this conjecture for  $sl_2$, see \fullref{conj thm}.

\subsection{The case of $sl_2$}\label{sl2 proof}
Let $L_1,\dots,L_n$ be irreducible finite-dimensional $sl_2$ modules
of highest weights $\La_1,\dots,\La_n\in\Z_{\geq 0}$.  Let
$v_1,\dots,v_n$ be the corresponding highest weight vectors.  Let
$V=L_1\otimes\dots\otimes L_n$. We also fix an $n$--tuple of non-zero
distinct complex numbers $\bs z=(z_1,\dots,z_n)$ and $l\in\Z_{\geq
  0}$. We set $\La_\infty=\sum_{s=1}^n\La_s-2l$.

In the case of $sl_2$ the weight function $\om(\bs t)$ can be rewritten in
the following form.
We say $\bs m=(m_1,\dots,m_n)\in \mc C(\bs \La,\La_\infty)$ if 
$m_s\in\{0,\dots,\La_i\}$, $s=1,\dots,n$, and $\sum_{s=1}^nm_s=l$. Set
\be
\om_{\bs m}(\bs t)=(\prod_{j=1}^n (m_j!)^{-1}) Sym \ \prod_{s=1}^n\ \ 
\prod_{i=m_1+\dots+m_{s-1}+1}^{m_1+\dots+m_s}\frac{1}{t_i-z_s},
\ee
where $Sym$ denotes the symmetrization with respect to $t_1,\dots,t_l$. Let 
\be
F^{\bs m}\bs v:=F^{m_1}v_1\otimes\dots \otimes F^{m_n}v_n.
\ee
Then we explicitly have 
\be
\om(\bs t)=\sum_{\bs m \in \mc C(\bs \La,\La_\infty)}
\om_{\bs m}(\bs t)F^{\bs m}\bs v.
\ee
Recall that if $\bs t$ is a critical point of the master function 
\eqref{master} then the vector $\om(\bs t)$ is called the Bethe vector.

It follows from \cite{ScV}, that there exists a Zariski open set
$U_1=U(\bs \La)$ in $\C^n$ such that for any $\bs z\in U_1$ there
exists a Zariski open set $U_2=U_2(\bs \La,\bs z)$ in $\C$ such for
all $\la\in U_2(\bs z)$, the number of orbits of critical points of
the $sl_2$ master function \eqref{master} associated to
$(\La_\infty,\la)$ equals to the dimension of the subspace of
$V$ of vectors of weight $\La_\infty=\sum_{s=1}^n\La_s-2l$.
Moreover all critical points are non-degenerate and the corresponding
Bethe vectors form a basis in this subspace.

\begin{theorem}\label{conj thm}
  Let $w$ be the generator of the $sl_2$ Weyl group.  
  For $\bs z\in U_1$, there exists a Zariski open set $U_3(\bs
  z)\subset\C$ with the following property. Let $\la\in U_3(\bs z)$, let
  $\bs t$ be a critical point associated to $(\La_\infty,\la)$ and 
  let $y$ be the corresponding polynomial. 
  Let $\bs t_w$ be the tuple represented by the polynomial $wy$. 
 
  Then all roots of the polynomial $wy$ are simple and  the vector
  $\mc A_{w}(\la)\,\om(\bs t)$ is a non-zero scalar multiple of the Bethe
  vector $\om(\bs t_w)$.
 \end{theorem}
\begin{proof}
  In the $sl_2$ case $\nu,\mu\in\C$, and the Casimir operator is given
  by $C=h\otimes h/2+e\otimes f+ f\otimes e$.
  
  We claim that the joint spectrum of $H_k(\la+\rho+\La_\infty/2)$, 
  $k=1,\dots,n$, acting in $V[\La_\infty]$, is
  generically simple. Indeed, in 
 the limit $\la \to \infty$ the main term is given by the
operators $\la h^{(k)}$.  The joint spectrum of commuting operators
$h^{(k)}$, $k=1,\dots,n$, is simple. Therefore $H_k(\la+\rho+\La_\infty/2)$, 
$k=1,\dots,n$,
for generic $\la$ have a simple joint spectrum as well.

It follows that the dynamical Weyl group maps the Bethe vectors to the
Bethe vectors.

Now we compare the two actions. We do it in the same limit $\la\to \infty$. 

The common eigenvectors of operators $h^{(i)}$ are monomial vectors
$F^{\bs m}\bs v$. The Weyl group of $sl_2$ is generated by the element
$w$, $w^2=id$ which acts on the weight vectors by
\bean\label{st act}
w (F^{m_1}v_1\otimes\dots \otimes F^{m_n}v_n) =
c F^{\La_1-m_1}v_1\otimes\dots\otimes
F^{\La_n-m_n}v_n,
\eean
where $c$ is some non-zero constant depending on $m_i,\La_i$.

By \fullref{proper}, the limit $\la\to \infty$, 
the dynamical Weyl group action on the Bethe vectors
coincides (up to a scalar) with the action of the Weyl group
\eqref{st act}.

Let us consider the limit of the action defined in terms of the
reproduction procedure.  It is shown in the proof of Theorem \eqref{sl2
  rep} that if $y$ represents an $sl_2$ critical point and $\tilde y$
is the immediate descendent, then for almost all $\la$, $y,\tilde y$
can be included in a family of critical points $y_a$, and their
descendents $\tilde y_a$ and in the limit $\la_a\to \infty$ the
product $y_a\tilde y_a$ tends to $T=\prod_{i=1}^n(t-z_i)^{\La_i}$.

Finally we claim that if the polynomials $y_a$ of degree
$l=\sum_{i=1}^n m_i$, represent critical points associated to $\la_a$
and the limit of $y_a$ as $\la_a$ tend to $\infty$ is
$\prod_i(x-z_i)^{m_i}$, then the corresponding Bethe vectors tend to a
scalar multiple of the monomial vector $F^{\bs m}\bs v$.

For $i=1,\dots,l$, let $s(i)\in \{1,\dots,n\}$ be such that the $i$-th
root of $y$, $t_i$, tends to $z_{s(i)}$. Then we write
$t_i(\la)=z_{s(i)}+c_i/\la+o(1/\la)$. The Bethe Ansatz equation for
$t_i$ implies that for any $j=1,\dots, n$, the set of $\{ c_i \ |\ 
s(i)=j\}$ satisfy the Bethe Ansatz equation with $n=1$: \be
-\frac{\La_j}{c_i}+\sum_{k,\ k\neq i,\ s(k)=s(i)}\frac{2}{c_i-c_k}=1.
\ee These equations are solved explicitly. The solutions are limits of
\cite[formulas (1.3.2)]{V} as $\beta\to \infty$. It follows that all
$c_i$ with $s(i)=j$ are different from zero and from each other.

Now consider the limit of the corresponding Bethe vector. 
The dominant term is 
\be
\la^l\prod_{j=1}^n (m_j!)^{-1}\prod_{i=1}^l c_i^{-1}F^{\bs m}\bs v.
\ee

This finishes the proof of the claim and the theorem.
\end{proof}

\section{Exponential populations}\label{exp sec}
We considered in detail the trigonometric Gaudin model, where the
Bethe Ansatz equation takes the form \eqref{Bethe eqn}. There are other
models, where the reproduction procedure for the solutions of the
Bethe Ansatz equation works in the same way and one obtains a
transitive and free Weyl group action on each population.  One such
model, the quasi-periodic Gaudin model, is considered in this section,
another one, the quasi-periodic $\XXX$ model, is considered in \fullref{xxx sec}.

We fix our $\g,\bs \La,\La_\infty,\bs l,\la$ as in \fullref{master
  sec}.  Let $z_1,\dots,z_n$ be any distinct complex numbers. Consider
the {\it master function with exponents}
\begin{gather}\label{master 2}
\Phi^{exp}(\bs t; \La_\infty;\la) = 
\prod_{i=1}^r\prod_{j=1}^{l_i}e^{-(\la,\al_i)t_j^{(i)}}
\prod_{i=1}^r\prod_{j=1}^{l_i}\prod_{s=1}^n
  (t_j^{(i)}-z_s)^{-(\La_s,\al_i)}\ \times
\\  
\prod_{i=1}^r\prod_{1\leq j<s\leq l_i} (t_j^{(i)}-t_s^{(i)})^{(\al_i,\al_i)}
  \prod_{1\leq i<j\leq r}
\prod_{s=1}^{l_i}\prod_{k=1}^{l_{j}}(t_s^{(i)}-t_k^{(j)})^{(\al_i,\al_j)}.
\notag
\end{gather}
We call $\bs t=(t_j^{(i)})^{j=1\dots l_i}_{i=1,\dots,r}$ {\it a
  critical point of the master function with exponents associated to
  $(\La_\infty,\la)$} if
\bean\label{Bethe eqn 2}
-(\la,\al_i)-\sum_{s=1}^n \frac{(\Lambda_s,
  \alpha_i)}{t_j^{(i)}-z_s}\ +\ 
\sum_{s,\ s\neq i}\sum_{k=1}^{l_s} \frac{(\alpha_s, \alpha_i)}{
  t_j^{(i)} -t_k^{(s)}}\ +\ 
\sum_{s,\ s\neq j}\frac {(\alpha_i, \alpha_i)}{ t_j^{(i)} -t_s^{(i)}}
= 0,
\eean
for $i=1,\dots,r$, $j=1,\dots,l_i$.

We have analogs of Propositions \ref{sln number} and \ref{sl2 number}.
\begin{prop}\label{sln number exp}
  Let $\g=sl_{N+1}$.  For almost all $\la$ the number of orbits of
  critical points of the master function with exponents associated to
  $(\La_\infty,\la)$ and counted with multiplicities does not exceed
  the dimension of the subspace of the weight $\La_\infty$ in the
  tensor product $L_{\La_1}\otimes\dots\otimes L_{\La_n}$.
\end{prop}
\begin{proof}
  The number of critical points of the master function with exponents
  \eqref{master 2} is finite for almost all $\la$, see \cite{MTV1}.
  
  Replacing the factors $e^{-(\la,\al_i)t_j^{(i)}}$ in the master
  function \eqref{master 2} with $(1+t_j^{(i)}/m)^{-(\la,\al_i)m}$ we
  obtain a master function of type \eqref{master}. Therefore the
  function \eqref{master 2} is the limit of master functions of type
  \eqref{master} as $m\to\infty$.  The proposition now follows from
  \fullref{sln number} and the fact that the number of orbits
  of isolated critical points of a function counted with multiplicity
  does not change under small deformations of the function.
\end{proof}

A different proof of \fullref{sln number exp} which uses Schubert
Calculus is given in \cite{MTV1}.

\begin{prop}\label{sl2 number exp}
  Let $\g=sl_2$.  For almost all $\la\in\C$ and almost all
  $(z_1,\dots,z_n)\in\C^n$, the number of orbits of critical points of
  the master function with exponents associated to $(\La_\infty,\la)$
  equals the dimension of the subspace of weight ${\La_\infty}$ in the
  tensor product $L_{\La_1}\otimes\dots\otimes L_{\La_n}$. Moreover
  all these points are non-degenerate.
\end{prop}
\begin{proof}
  If $\la$ is a large positive integer then the proposition is proved
  by using methods of \cite{MV6}.  The rest is similar to the proof of
  \fullref{sl2 number}
\end{proof}

Let $\g$ be a Kac--Moody algebra.  As in \fullref{gen pop sec} we
represent a tuple $\bs t=(t_j^{(i)})_{i=1,\dots,r}^{j=1,\dots,l_i}$ by
the $r$--tuple of polynomials $\bs y=(y_1,\dots,y_r)$, where
$y_i=\prod_{j=1}^{l_i}(x-t_j^{(i)})$, $i=1,\dots,r$. We make no 
distinction between $(y_1,\dots,y_r)$ and $(c_1y_1,\dots,c_ry_r)$
where $c_1,\dots, c_r$ are non-zero complex numbers.  We introduce
polynomials $T_i$, $i=1,\dots,r$, by formula \eqref{T}.

We call an $r$--tuple of polynomials $\bs y$ {\it exponentially
  off-diagonal} if its roots do not belong to the union of singular
hyperplanes in \eqref{Bethe eqn 2}.  Namely $\bs y$ is exponentially
off-diagonal if for $i=1,\dots,r$, all roots of the polynomial $y_i$
are simple, different from the roots of the polynomials $y_j$ for all
$j$ such that $(\al_j,\al_i)\neq 0$ and different from the roots of
the polynomial $T_i$.

We have the corresponding exponential reproduction procedure. Namely,
an $r$--tuple of polynomials $\bs y$ is called {\it exponentially
  fertile in the $i$-th direction with respect to $\la$},
$i\in\{1,\dots,r\}$, if there exists a polynomial $\tilde y_i$ such
that
\be
W(y_i,e^{\langle \la,\al_i^\vee\rangle x}\tilde y_i)=e^{\langle
  \la,\al_i^\vee\rangle x}T_i\prod_{j=1, j\neq i}^r y_j^{-a_{ij}}.
\ee
Then the $r$--tuple of polynomials $y^{(i)}=(y_1,\dots,\tilde
y_i,\dots,y_r)$ is called {\it an exponential immediate descendent of
  $\bs y$ with respect to $\la$ in the direction $i$}.

An $r$--tuple of polynomials is called {\it exponentially fertile with
  respect to $\la$} if it is exponentially fertile with respect to
$\la$ in all directions $i=1,\dots,r$.

\begin{theorem} Let $\g$ be a Kac--Moody algebra.
  An $r$--tuple of polynomials $\bs y$ represents a $\g$ critical point
  of the master function with exponents associated to
  $(\La_\infty,\la)$ if and only if $\deg y_i=l_i$, $\bs y$ is
  exponentially off-diagonal and exponentially fertile with respect to
  $\la$.  Moreover, if $\bs y$ represents a critical point of the
  master function with exponents associated to $(\La_\infty,\la)$ and
  if the immediate descendent of $\bs y$ with respect $\la$ in the
  $i$-th direction, $y^{(i)}=(y_1,\dots,\tilde y_i,\dots,y_r)$, is
  exponentially off-diagonal then $y^{(i)}$ represents a critical
  point of the master function with exponents associated to
  $(s_i\La_\infty,s_i\la)$.
\end{theorem} 
\begin{proof}
  The proof is similar to the proof of \fullref{gen rep}.
\end{proof}

An $r$--tuple of polynomials $\bs y^{(i_1,i_2,\dots,i_m)}$, where
$i_k\in\{1,\dots,r\}$, $k=1,\dots,m$, is called {\it an exponential
  descendent of $\bs y$ with respect to $\la$ in the directions
  $(i_1,\dots,i_m)$} if there exist $r$--tuples of polynomials $\bs
y^{(i_1,i_2,...,i_k)}$, $k=1,\dots,m-1$, such that for $k=1,\dots,m$,
the $r$--tuple $\bs y^{(i_1,i_2,...,i_k)}$ is an exponential immediate
descendent of $\bs y^{(i_1,i_2,...,i_{k-1})}$ with respect to
$s_{i_{k-1}}\dots s_{i_2} s_{i_1} \la$ in the $i_k$-th direction.

An $r$--tuple of polynomials $\bs y$ is called {\it exponentially
  super-fertile with respect to $\la$} if it is exponentially fertile
with respect to $\la$ and all exponential descendents $\bs
y^{(i_1,i_2,\dots,i_m)}$ of $\bs y$ with respect to $\la$ in the
directions $(i_1,\dots,i_m))$ are exponentially fertile with respect
to $s_{i_m}\dots s_{i_1}\la$.

For any $N$--tuple of functions $\bs y$ and an $sl_{N+1}$ weight $\la$,
we set $y_{N+1}=1$ and define the linear differential operator of
order $N+1$:
\be
D^{exp}(\bs y,\la)= \prod^{N\ \to\ 0}_i\left(
\partial \ -\ \ln' \ \left( 
\frac{ y_{i+1}\prod_{s=1}^{i} e^{(\la,\al_s) x}
  T_s }{ y_{i} }\ \right)\right).
\ee

\begin{prop}
  Let $\bs y$ be an $N$--tuple of polynomials and $l_i=\deg y_i$, 
  $i=1,\dots,N$. Let $\bs \La_\infty$ be given by \eqref{at inf}. 
  Let $\bs y$ represent an $sl_{N+1}$ critical point of the master
  function with exponents associated to $(\La_\infty,\la)$ or let $\bs
  y$ be exponentially super-fertile with respect to $\la$.
  Then the kernel of the operator $D^{exp}(\bs y,\la)$ is spanned 
  by functions of the form
\bean\label{exp kernel form}
p_0, p_1e^{(\la,\al_1)x}, \dots, p_N e^{(\la,\al_1+\dots+\al_N)x},
\eean
where $p_i$ is a polynomial of degree $\deg
y_1+(\La_\infty,\al_1+\dots+\al_i)$.  The only singular points of the
operator $D^{exp}(\bs y,\la)$ in $\C$ are regular singular points
located at $z_1,\dots,z_n$, and the exponents at $z_i$, $i=1,\dots,n$,
are
\bean\label{exp exponents}
z_i:  0, 
(\La_i+\rho,\al_1), 
(\La_i+\rho,\al_1+\al_2),\dots,(\La_i+\rho,\al_1+\dots+\al_N),
\eean

Conversely, if a linear differential operator $D$ of order $N+1$ has
the kernel spanned by functions of the form \eqref{exp kernel form} and
the only singular points of $D$ in $\C$ are regular singular points at
$z_1,\dots,z_n$ with the exponents given by \eqref{exp exponents}, then
the $N$--tuple $\bs y$ given by the divided Wronskians
\be
y_i=
\frac{W(p_0, p_1e^{(\la,\al_1)x}, \dots, 
p_{i-1}e^{(\la,\sum_{j=1}^{i-1}\al_j)x})}
{e^{(\la,\sum_{j=1}^{i-1}(i-j)\al_{j})x}\prod_{j=1}^{i-1}T_j^{i-j}},
\ee
$i=1,\dots,N$, is an $N$--tuple of polynomials  
which is exponentially super-fertile with respect to $\la$ and satisfies $\deg y_i=l_i$, $i=1,\dots,N$.
\end{prop}
\begin{proof} The proof is similar to the proof of
  \fullref{sing points}.
\end{proof}

\begin{conj}\label{exp conj} If an $r$--tuple of polynomials 
  $\bs y$ represents a critical point of the master function with
  exponents associated to $(\La_\infty,\la)$ then $\bs y$ is
  exponentially super-fertile with respect to $\la$.
\end{conj}

\begin{theorem}
  \fullref{exp conj} holds for the case of simple Lie algebras
  of types $A_N$ and $B_N$.
\end{theorem}
\begin{proof}
  The proof is similar to the proof of Theorems \ref{super},
  \ref{super B}.
\end{proof}

For an $r$--tuple of polynomials $\bs y$ and a $\g$ weight $\la$, we
denote $P^{exp}(\bs y,\la)$ the set of all pairs of the form $(\bs
y^{(i_1,i_2,\dots,i_m)}, s_{i_m} \dots s_{i_2}s_{i_1} \la)$, where
$m\in\Z_{\geq 0}$, $i_k\in\{1,\dots,r\}$, $k=1,\dots,m$, and $\bs
y^{(i_1,i_2,...,i_m)}$ is an exponential descendent of $\bs y$ with
respect to $\la$ in the directions $(i_1,\dots,i_m)$.

We call the set $P^{exp}(\bs y,\la)$ {\it the exponential
  prepopulation originated at $(\bs y,\la)$}.  Let an $r$--tuple of
polynomials $\bs y$ be exponentially super-fertile with respect to
$\la$.  We call the set $P^{exp}(\bs y,\la)$ {\it the exponential
  population originated at $(\bs y,\la)$}.

\begin{theorem}
  Let $\g$ be any simple Lie algebra and let $\la$ be a strongly
  non-integral $\g$--weight. Let an $r$--tuple of polynomials $\bs y$ be
  exponentially super-fertile with respect to $\la$. Then the map
  $P^{exp}(\bs y,\la)\to \mc W\la$ such that $(\tilde {\bs y},\tilde
  \la)\mapsto \tilde \la$ is a bijection of the exponential population
  originated at $(\bs y,\la)$ and of the orbit of the Weyl group.
\end{theorem}
\begin{proof}
  The proof is similar to the proof of \fullref{general cor}.
\end{proof}

\section{Difference reproduction}\label{xxx sec}
In this section we describe the Bethe Ansatz equation corresponding to
the quasi-periodic $\XXX$ model. In this case the reproduction
procedure works in a similar way and one obtains a free and transitive
Weyl group action on a population.

Let $h\in \C$ be a complex non-zero number.  We fix $\g,\bs
\La,\La_\infty,\bs l,\la$ as in \fullref{master sec}.  Let
$z_1,\dots,z_n$ be any distinct complex numbers, subject to the
conditions $z_i-z_j\not\in h\Z$ for all $i,j\in\{1,\dots,n\}$, $i\neq
j$.

Consider the {\it exponential $\XXX$ Bethe equation} on
variables $\bs t=(t_j^{(i)})_{i=1,\dots,r}^{j=1,\dots, l_i}$:
\begin{align}\label{XXX BAE}
e^{\langle\la,\al_i^\vee\rangle h}=&
\prod_{s=1}^n  \frac{t_j^{(i)}-z_s+(\La_s,\al_i) h/2}
{t_j^{(i)}-z_s-(\La_s,\al_i) h/2}\ \times\\
&\prod_{m=1,\dots,r,\atop m\neq i }\hspace{-7pt}\left(\prod_{k=1}^{l_m}
\frac{t_j^{(i)}-t_k^{(m)}+ h/2}
{t_j^{(i)}-t_k^{(m)}-h/2}\right)^{-a_{im}}\hspace{-10pt} 
\prod_{k=1,\dots,l_i, \atop k\neq j}
\frac{t_j^{(i)}-t_k^{(i)}-h}{t_j^{(i)}-t_k^{(i)}+h},\notag
\end{align}
where $i=1,\dots,r$, $j=1,\dots, l_i$.
 
As in \fullref{gen pop sec} we represent a tuple $\bs
t=(t_j^{(i)})_{i=1,\dots,r}^{j=1,\dots,l_i}$ by the $r$--tuple of
polynomials $\bs y=(y_1,\dots,y_r)$, where
$y_i=\prod_{j=1}^{l_i}(x-t_j^{(i)})$, $i=1,\dots,r$. We make no 
distinction between $(y_1,\dots,y_r)$ and $(c_1y_1,\dots,c_ry_r)$
where $c_1,\dots, c_r$ are non-zero complex numbers.

For $i=1,\dots,r$, set
\be
T_i^{(h)}(x)=\prod_{s=1}^n\prod_{j=1}^{(\La_s,\al_i)}
(x-z_s-(\La_s,\al_i) h/2+jh).
\ee
An $r$--tuple of polynomials $\bs y$ is called {\it exponentially
  difference off-diagonal with respect to $(\bs \La, \bs z;h)$} if for
$i=1,\dots,r$ the polynomial $y_i(x)$ has only simple roots, different
from the roots of polynomials $y_m(x+h/2)$, whenever
$(\al_i,\al_m)\neq 0$, and different from the roots of polynomials
$T_i^{(h)}$, $y_i(x+h)$.

A solution $\bs t$ of \eqref{XXX BAE} is called {\it off-diagonal} if
the corresponding $r$--tuple of polynomials is exponentially difference
off-diagonal.

\begin{lem}  A polynomial $y$ of degree $l$ represents an $sl_2$
 off-diagonal solution of exponential $\XXX$ Bethe equation associated to
  $(\La_\infty,\la)$  
  if and only if $y$ is exponentially difference off-diagonal
  and there exists a polynomial $B(x)$ such that
\begin{align*}
y(x+h)e^{\langle\la,\al^\vee\rangle h}&\prod_{s=1}^n (x-z_s-\frac{(\La_s,\al) h}{2}) 
\\
&+ B(x)y(x)+ y(x-h)\prod_{s=1}^n(x-z_s+\frac{(\La_s,\al) h}{2}) =0. 
\end{align*}\end{lem}
\begin{proof}
The lemma is proved similarly to \fullref{st}.
\end{proof}

\begin{prop}\label{sln number xxx} Let $\g=sl_{N+1}$.
  For almost all $\la$ the number of orbits of off-diagonal solutions
  of the exponential $\XXX$ Bethe Ansatz equations associated to
  $(\La_\infty,\la)$ does not exceed the dimension of the subspace of
  the weight $\La_\infty$ in the tensor product
  $L_{\La_1}\otimes\dots\otimes L_{\La_n}$.
\end{prop}
\begin{proof}
  The proof is similar to the proof of the \fullref{sln number
    exp} with the help of \cite[Corollary 4.15]{MV3}.
\end{proof}

\begin{prop}\label{sl2 number xxx}{\rm\cite{TV}}\qua Let $\g=sl_2$.
  For almost all $\la$ and almost all $(z_1,\dots,z_n)\in\C^n$, the
  number of orbits of solutions $\bs t$ of the exponential $\XXX$ Bethe
  Ansatz equation associated to $(\La_\infty,\la)$ such that $t_i\neq
  t_j$ equals the dimension of the subspace of weight ${\La_\infty}$
  in the tensor product $L_{\La_1}\otimes\dots\otimes L_{\La_n}$.
  Moreover all such solutions are non-degenerate.
\end{prop}

We now describe the corresponding exponential difference reproduction
procedure.

Denote $W_h$ the discrete Wronskian:
\be
W_h(f_1,\dots,f_N):=\det({f_i(x+(j-1)h)})_{i,j=1,\dots,N}.
\ee
An $r$--tuple of polynomials $\bs y$ is called {\it exponentially
  difference fertile with respect to $\la$ in the $i$-th direction},
$i\in\{1,\dots,r\}$, if there exists a polynomial $\tilde y_i$ such
that
\be
W_{h}(y_i,e^{\langle\la,\al_i^\vee\rangle x}\tilde y_i)=
e^{\langle\la,\al_i^\vee\rangle x}\ T_i^{(h)}(x)\prod_{m=1,\ m\neq
  i}^r (y_m(x+h/2))^{-a_{im}}.  
\ee 
Then the $r$--tuple of polynomials $y^{(i)}=(y_1,\dots,\tilde
y_i,\dots,y_r)$ is called {\it an exponential difference immediate
  descendent of $\bs y$ with respect to $\la$ in the $i$-th
  direction}.

An $r$--tuple is called {\it exponentially difference fertile with
  respect to $\la$} if it is exponentially difference fertile with
respect to $\la$ in all directions $i=1,\dots,r$.

\begin{theorem}
  An $r$--tuple of polynomials $\bs y$ represents an off-diagonal
  solution of the exponential $\XXX$ Bethe Ansatz equation associated
  to $(\La_\infty,\la)$ if and only if $\bs y$ is exponentially
  difference off-diagonal, $\deg y_i=l_i$, $i=1,\dots,r$, and $\bs y$
  is exponentially difference fertile with respect to $\la$.
  Moreover, if $\bs y$ represents an off-diagonal solution of the
  exponential $\XXX$ Bethe Ansatz equation associated to
  $(\La_\infty,\la)$ and if the exponential difference immediate
  descendent of $\bs y$ with respect $\la$ in the $i$-th direction,
  $y^{(i)}=(y_1,\dots,\tilde y_i,\dots,y_r)$, is exponentially
  difference off-diagonal then $y^{(i)}$ represents an off-diagonal
  solution of the exponential $\XXX$ Bethe Ansatz equation associated
  to $(s_i\La_\infty,s_i\la)$.
\end{theorem} 
\begin{proof}
  The proof is similar to the proof of \fullref{gen rep}, cf
  \cite{MV5} and also \cite{MV2,MV3}.
\end{proof}

An $r$--tuple of polynomials $\bs y^{(i_1,i_2,\dots,i_m)}$, where
$m\in\Z_{\geq 0}$, $i_k\in\{1,\dots,r\}$, $k=1,\dots,m$, is called
{\it an exponential difference descendent of $\bs y$ with respect to
  $\la$ in the directions $(i_1,\dots,i_m)$} if there exist $r$--tuples
of polynomials $\bs y^{(i_1,i_2,...,i_k)}$, $k=1,\dots,m-1$, such that
for $k=1,\dots,m$, the $r$--tuple $\bs y^{(i_1,i_2,...,i_k)}$ is an
exponential difference immediate descendent of $\bs
y^{(i_1,i_2,...,i_{k-1})}$ with respect to $s_{i_{k-1}}\dots s_{i_2}
s_{i_1} \la$ in the $i_k$-th direction.

An $r$--tuple of polynomials $\bs y$ is called {\it exponentially
  difference super-fertile with respect to $\la$} if it is
exponentially difference fertile with respect to $\la$ and all
exponential descendents $\bs y^{(i_1,i_2,\dots,i_m)}$ of $\bs y$ with
respect to $\la$ in the directions $(i_1,i_2,\dots,i_m)$ are
exponentially difference fertile with respect to $s_{i_m}\dots
s_{i_1}\la$.

For any $N$--tuple of functions $\bs y$ and an $sl_{N+1}$ weight $\la$,
we set $y_{N+1}=1$ and define the linear difference operator:
\begin{align*}
D_h^{exp}(\bs y,\la)=&
\prod^{N \to 0}_i
\Bigg( \partial_h  -  
\frac{y_{i+1}(x+(i+2)h/2)}{y_{i+1}(x+ih/2)} 
\frac{y_{i}(x+(i-1)h/2)}{y_{i}(x+(i+1)h/2)}\ \times\\&
\prod_{s=1}^{i} \frac{e^{h(\la,\al_s)}T_{s}(x + (2i-s+1)h/2)}
{T_{s}(x+(2i-s-1)h/2)}
 \Bigg),
\end{align*}
where $\partial_h$ is the operator acting on functions of $x$ by the
formula $\partial_h(f(x))=f(x+h)$.

Let $V$ be a space spanned by functions of the type $p_0e^{\la_0 x},
p_1e^{\la_1x}, \dots, p_N e^{\la_N x}$ where $p_i$, $i=0,\dots,N$, are
polynomials and $\la_i\in\C$, $i=0,\dots,N$.  We say the space $V$
{\it has no base points} if for any $z\in\C$ there exists $f\in V$,
such that $f(z)\neq 0$.

Assume $V$ has no base points. For $i=2,\dots,N$, let $U_i$ be the
monic polynomial of the greatest possible degree such that
$W_h(f_1,\dots, f_i)/U_i$ is a holomorphic function for all
$f_1,\dots,f_i\in V$.  Following \cite{MV3}, we call an $N$--tuple of
monic polynomials $(T_1,\dots,T_N)$ {\it a frame of space $V$} if for
$i=2,\dots,N$ we have $U_i=\prod_{j=1}^{i-1}\prod_{s=1}^{i-j}T_j(x+(s-1)h)$.

\begin{lem}
  Let $V$ be a space spanned by functions of the type $$p_0e^{\la_0 x},
  p_1e^{\la_1x}, \dots, p_N e^{\la_N x}.$$ Let $V$ have no base
  points.  Then there exists a unique frame of $V$.
\end{lem}
\begin{proof}
  The proof is similar to the proof of \cite[Lemma 4.9]{MV3}.
\end{proof}

\begin{prop}
    Let $\bs y$ be an $N$--tuple of polynomials and $l_i=\deg y_i$, 
  $i=1,\dots,N$. Let $\bs \La_\infty$ be given by \eqref{at inf}. 
   Let $\bs y$ represent an off-diagonal solution of $sl_{N+1}$
  exponential $\XXX$ Bethe Ansatz equation associated to
  $(\La_\infty,\la)$ or let $\bs y$ be exponentially difference
  super-fertile with respect to $\la$. Then the kernel of the
  operator $D_h^{exp}(\bs y,\la)$ is spanned by functions of the form
\bean\label{exp diff kernel form}
p_0, p_1e^{(\la,\al_1)x}, \dots, p_N e^{(\la,\al_1+\dots+\al_N)x},
\eean
where $p_i$ is a polynomial of degree 
$\deg y_1+(\La_\infty,\al_1+\dots+\al_i)$. Moreover, the $N$--tuple
\bean\label{frame}
(T_1^{(h)}(x),T_2^{(h)}(x+h/2),\dots,T_N^{(h)}(x+(N-1)h/2)) 
\eean
is the frame of the kernel of the operator $D_h^{exp}(\bs y,\la)$.

Conversely, if a linear difference operator $D$ of order $N+1$ has the
kernel spanned by functions of the form \eqref{exp diff kernel form}
with the frame \eqref{frame} then the $N$--tuple $\bs y$ given by
\be
y_i=
\frac{W_h(p_0, p_1e^{(\la,\al_1)x}, \dots, 
p_{i-1}e^{(\la,\sum_{j=1}^{i-1}\al_j)x})}
{e^{(\la,\sum_{j=1}^{i-1}(i-j)\al_{j})x}
\prod_{s=1}^{i-j}T_j^{(h)}(x+(s+j/2-3/2)h)},
\ee
$i=1,\dots,N$, is an $N$--tuple of polynomials which is exponentially
difference super-fertile with respect to $\la$ and satisfies $\deg
y_i=l_i$, $i=1,\dots,N$.
\end{prop}
\begin{proof}
  The proof is similar to the proof of \fullref{sing points}.
\end{proof}

\begin{conj}\label{exp diff conj} Let $\g$ be any simple Lie algebra.
  If an $r$--tuple of polynomials $\bs y$ represents an off-diagonal
  solution of the exponential $\XXX$ Bethe Ansatz equation associated
  to $(\La_\infty,\la)$ then $\bs y$ is exponentially difference
  super-fertile with respect to $\la$.
\end{conj}

\begin{theorem}
  \fullref{exp diff conj} holds for the case of simple Lie
  algebras of types $A_N$ and $B_N$.
\end{theorem}
\begin{proof}
  The proof is similar to the proof of Theorems \ref{super},
  \ref{super B}.
\end{proof}

For an $r$--tuple of polynomials $\bs y$ we denote $P^{exp}_h(\bs
y,\la)$ the set of all pairs of the form \newline $(\bs
y^{(i_1,i_2,\dots,i_m)}, s_{i_m} \dots s_{i_2}s_{i_1} \la)$, where
$m\in\Z_{\geq 0}$, $i_k\in\{1,\dots,r\}$, $k=1,\dots,m$, and $\bs
y^{(i_1,i_2,...,i_m)}$ is an exponential difference descendent of $\bs
y$ with respect to $\la$ in directions $(i_1,\dots,i_m)$.

We call the set $P^{exp}_h(\bs y,\la)$ {\it the exponential difference
  prepopulation originated at $(\bs y,\la)$}.  If an $r$--tuple of
polynomials $\bs y$ is exponentially difference super-fertile with
respect to $\la$, then we call the set $P_h^{exp}(\bs y,\la)$ {\it the
  exponential difference population originated at $(\bs y,\la)$}.

\begin{theorem}
  Let $\g$ be any simple Lie algebra and let $\la$ be a strongly
  non-integral $\g$--weight. Let an $r$--tuple of polynomials $\bs y$ be
  exponentially difference super-fertile with respect to $\la$. Then
  the map $P_h^{exp}(\bs y,\la)\to \mc W\la$ such that $(\tilde {\bs
    y},\tilde \la)\mapsto \tilde \la$ is a bijection of the exponential
  difference population originated at $(\bs y,\la)$ and of the orbit
  of the Weyl group.
\end{theorem}
\begin{proof}
  The proof is similar to the proof of \fullref{general cor}.
\end{proof}

\bibliographystyle{gtart}
\bibliography{link}

\end{document}